\documentclass[12pt]{article}

\usepackage{amsmath, amssymb, ntheorem,booktabs, color}
\usepackage{algorithm}
\usepackage{algorithmic}
\usepackage{rotating}
\usepackage{tabulary}
\usepackage{bm}
\usepackage{caption}
\usepackage{url}
\usepackage{latexsym}
\def\qed{\hfill $\Box$}

\usepackage{subfig}

\usepackage{subfloat}

\newtheorem{thm}{Theorem}[section]
\newtheorem{prop}{Proposition}[section]
\newtheorem{lem}{Lemma}[section]

\numberwithin{equation}{section}

\begin{document}
\makeatletter

\begin{center}
\large{\bf Appropriate Learning Rates of Adaptive Learning Rate Optimization Algorithms for Training Deep Neural Networks}\\
{\small This work was supported by JSPS KAKENHI Grant Number JP18K11184.}
\end{center}\vspace{3mm}
\begin{center}

\textsc{Hideaki Iiduka}\\
Department of Computer Science, Meiji University
              1-1-1 Higashimita, Tama-ku, Kawasaki-shi, Kanagawa, 214-8571 Japan.
              (iiduka@cs.meiji.ac.jp)
\end{center}

\vspace{2mm}

\footnotesize{
\noindent\begin{minipage}{14cm}
{\bf Abstract:}
This paper deals with nonconvex stochastic optimization problems in deep learning  and provides appropriate learning rates with which adaptive learning rate optimization algorithms, such as Adam and AMSGrad, can approximate a stationary point of the problem. In particular, constant and diminishing learning rates are provided to approximate a stationary point of the problem. Our results also guarantee that the adaptive learning rate optimization algorithms can approximate global minimizers of convex stochastic optimization problems.  The adaptive learning rate optimization algorithms are examined in numerical experiments on text and image classification. The experiments show that the algorithms with constant learning rates perform better than ones with diminishing learning rates. 
\end{minipage}
 \\[5mm]


\hbox to14cm{\hrulefill}\par

\section{Introduction}
\label{sec:1}
The main objective of the field of deep learning is to train deep neural networks \cite{shao2014}, \cite{pas2019}, \cite{wu2019}, \cite{zhao2019} appropriately. One way of achieving the objective is to devise useful methods for finding model parameters of deep neural networks that reduce certain cost functions called the expected risk and empirical risk (Section 2 in \cite{bottou}). Accordingly, optimization methods are needed for minimizing the expected (or empirical) risk, i.e., for solving stochastic optimization problems in deep learning. 
 
The classical method for solving a convex stochastic optimization problem is the stochastic approximation (SA) method \cite{robb1951}, \cite{nem2009} which is a first-order method using the stochastic (sub)gradient of an observed function at each iteration. Modifications of the SA method, such as the mirror descent SA method \cite{nem2009} and the accelerated SA method \cite{gha2012}, have been presented.

As the field of deep learning has developed, useful algorithms based on the SA method and incremental methods \cite{nedic2001} have been presented to adapt the {\em learning rates} of all model parameters. These algorithms are called {\em adaptive learning rate optimization algorithms} (Subchapter 8.5 in \cite{deep}). For example, some algorithms use momentum (Subchapter 8.3.2 in \cite{deep})  or Nesterov's accelerated gradients (Subchapter 2.2 in \cite{nesterov} and Subchapter 8.3.3 in \cite{deep}). The AdaGrad algorithm \cite{adagrad} is a modification of the mirror descent SA method, while the RMSProp algorithm (Algorithm 8.5 in \cite{deep}) is based on AdaGrad. AdaGrad and RMSProp both use element-wise squared values of the stochastic (sub)gradient.

The Adam algorithm \cite{adam}, which is based on momentum and RMSProp, is a powerful algorithm for training deep neural networks. The performance measure of adaptive learning rate optimization algorithms is called the regret (see \eqref{REG} for the definition of regret), and the main objective of adaptive learning rate optimization algorithms is to achieve low regret. However, there is an example of a convex optimization problem in which Adam does not minimize the regret (Theorems 1--3 in \cite{reddi2018}). 

The AMSGrad algorithm \cite{reddi2018} was presented to guarantee the regret is minimized and preserve the practical benefits of Adam. In particular, AMSGrad must use diminishing learning rates (Theorem 4 in \cite{reddi2018}) to optimize deep neural network models. When adaptive learning rate optimization algorithms with diminishing learning rates are applied to complicated stochastic optimizations, the learning rates are approximately zero for a number of iterations, which implies that using diminishing learning rates would not be implementable in practice. Even if algorithms with diminishing learning rates could be made to work, we would need to empirically select suitable learning rates to increase their convergence speed. However, it is too difficult to select in advance suitable diminishing learning rates that guarantee sufficiently quick convergence since the selection significantly affects the model parameters (see Subchapter 8.5 in \cite{deep}). 

Another issue of adaptive learning rate optimization algorithms is that they cannot be applied to nonconvex stochastic optimization problems, while they can minimize the regret and achieve a low regret only when the cost functions are convex. Since the primary goal of training deep models is to solve nonconvex stochastic optimization problems in deep learning by using optimization algorithms, we need to develop optimization algorithms that can be applied to nonconvex stochastic optimization.

In this paper, we propose an adaptive learning rate optimization algorithm (Algorithm \ref{algo:1}) that can be applied to nonconvex stochastic optimization in deep learning. The advantage of the proposed algorithm is that it uses {\em constant} learning rates. In the case of constant learning rates, we can show that it approximately finds a stationary point of the nonconvex stochastic optimization problem (Theorem \ref{thm:1}). We also discuss how to set appropriate constant learning rates to find the stationary point. We show that the proposed algorithm can be applied to nonconvex stochastic optimization problems from the viewpoints of both theory and practice.

The proposed algorithm can also use {\em diminishing} learning rates. We provide sufficient conditions for the diminishing learning rates to ensure that the algorithm can solve the nonconvex stochastic optimization problem (Theorem \ref{thm:2}). We also determine the rate of convergence of the algorithm with diminishing learning rates to establish its performance. 

This paper makes three contributions. The first contribution of this paper is to enable us to consider stationary point problems associated with nonconvex stochastic optimization problems in deep learning, in contrast to the previously reported results in \cite{adam} and \cite{reddi2018} that presented algorithms for convex optimization. This implies that the results in this paper can be applied to nonconvex stochastic optimization problems in convolutional neural networks (CNNs) and their variants such as the residual network (ResNet).

The second contribution is to propose an adaptive learning rate optimization algorithm (Algorithm \ref{algo:1}) for solving the problem, together with its convergence analysis for constant learning rates and diminishing learning rates. Since constant learning rates are not zero for a number of iterations, the analysis for constant learning rates would be useful from the viewpoints of both theory and practice. In the special case where cost functions are convex, our analyses guarantee that the proposed algorithm can solve the convex stochastic optimization problem (Propositions \ref{prop:0} and \ref{prop:1}), in contrast to the previously reported results in \cite{adam} and \cite{reddi2018} showing that Adam and AMSGrad achieve low regret. 

The third contribution of this paper is that we show that the proposed algorithm can be applied to stochastic optimization with image and text classification tasks. The numerical results show that the proposed algorithm with constant learning rates performs better than the one with diminishing learning rates (Section \ref{sec:5}). 

This paper is organized as follows. Section \ref{sec:2} gives the mathematical preliminaries and states the main problem. Table \ref{notation} summarizes the notation used in this paper. Section \ref{sec:3} presents the adaptive learning rate optimization algorithm for solving the main problem and analyzes its convergence. Section \ref{sec:5} numerically compares the behaviors of the proposed algorithm with constant learning rates and with diminishing learning rates. Section \ref{sec:6} concludes the paper with a brief summary.

\section{Stationary Point Problem Associated With Nonconvex Optimization Problem}
\label{sec:2}

\begin{table*}[htbp]
\begin{center}
\caption{List of Notation}
\begin{tabular}{c||l} \hline
Notation & Description\\ \hline
$\mathbb{N}$ & 
The set of all positive integers and zero \\
$\mathbb{R}^d$ & 
A $d$-dimensional Euclidean space with inner product $\langle \cdot, \cdot \rangle$,\\ 
& which induces the norm $\| \cdot \|$ \\
$\mathbb{S}^d$ & 
The set of $d \times d$ symmetric matrices, i.e.,\\ 
& $\mathbb{S}^d = \{ M \in \mathbb{R}^{d \times d} \colon M = M^\top \}$ \\
$\mathbb{S}_{++}^d$ &
The set of $d \times d$ symmetric positive-definite matrices, i.e.,\\
& $\mathbb{S}_{++}^d = \{ M \in \mathbb{S}^{d} \colon M \succ O \}$ \\
$\mathbb{D}^d$ & 
The set of $d \times d$ diagonal matrices, i.e.,\\ 
& $\mathbb{D}^d = \{ M \in \mathbb{R}^{d \times d} \colon M = \mathsf{diag}(x_i), \text{ } x_i \in \mathbb{R} \text{ } (i=1,2,\ldots,d) \}$ \\
$A \odot B$ & 
The Hadamard product of matrices $A$ and $B$\\ 
& ($\bm{x} \odot \bm{x} := (x_i^2) \in \mathbb{R}^d$ ($\bm{x} := (x_i) \in \mathbb{R}^d$))\\
$\langle \bm{x}, \bm{y} \rangle_H$ &
The $H$-inner product of $\mathbb{R}^d$, where $H \in \mathbb{S}_{++}^d$, 
i.e.,
$\langle \bm{x}, \bm{y} \rangle_H := \langle \bm{x}, H \bm{y} \rangle$ \\
$\|\bm{x}\|_H^2$ & 
The $H$-norm, where $H \in \mathbb{S}_{++}^d$, i.e., $\|\bm{x}\|_H^2 := \langle \bm{x}, H \bm{x} \rangle$ \\
$P_X$ &
The metric projection onto a nonempty, closed convex set $X$ $(\subset \mathbb{R}^d)$ \\
$P_{X,H}$ &
The metric projection onto $X$ under the $H$-norm \\
$\mathbb{E}[Y]$ & 
The expectation of a random variable $Y$ \\
$\bm{\xi}$ &
A random vector whose probability distribution $P$\\ 
& is supported on a set $\Xi \subset 
\mathbb{R}^{d_1}$ \\
$F(\cdot, \bm{\xi})$ &
A continuously differentiable function from $\mathbb{R}^d$ to $\mathbb{R}$\\
& for almost every 
$\bm{\xi} \in \Xi$ \\
$f$ & 
The objective function defined for all $\bm{x}\in \mathbb{R}^d$ by $f(\bm{x}) := \mathbb{E}[F(\bm{x}, \bm{\xi})]$ \\
$\nabla f$ &
The gradient of $f$ \\
$\mathsf{G}(\bm{x},\bm{\xi})$ & 
Stochastic gradient for a given $(\bm{x},\bm{\xi}) \in \mathbb{R}^d \times \Xi$\\
& which satisfies $\mathbb{E}[\mathsf{G}(\bm{x},\bm{\xi})]= \nabla f(\bm{x})$
\\
$X^\star$ & 
The set of stationary points of the problem of minimizing $f$ over $X$ \\ \hline
\end{tabular}\label{notation}
\end{center}
\end{table*}

Let us consider the following problem (see, e.g., Subchapter 1.3.1 in \cite{facc1} for the details of stationary point problems):

\textbf{Problem:} Assume that 
\begin{enumerate} 
\item[(A1)] $X \subset \mathbb{R}^d$ is a closed convex set onto which the projection can be easily computed; 
\item[(A2)] $f \colon \mathbb{R}^d \to \mathbb{R}$ defined for all $\bm{x}\in \mathbb{R}^d$ by $f(\bm{x}) := \mathbb{E}[F(\bm{x}, \bm{\xi})]$ is well defined, where $F(\cdot,\bm{\xi})$ is continuously differentiable for almost every $\bm{\xi} \in \Xi$.
\end{enumerate}
Then, find a  stationary point $\bm{x}^\star$ of the problem of minimizing $f$ over $X$, i.e.,
\begin{align*}
\bm{x}^\star \in X^\star := 
\left\{ \bm{x}^\star \in X \colon 
\langle \bm{x} - \bm{x}^\star, \nabla f (\bm{x}^\star) \rangle \geq 0
\text{ } (\bm{x}\in X)
\right\}.
\end{align*}

If $X = \mathbb{R}^d$, then $X^\star = \{ \bm{x}^\star \in \mathbb{R}^d \colon \nabla f(\bm{x}^\star) = \bm{0} \}$. If $f$ is convex, then $\bm{x}^\star \in X^\star$ is a global minimizer of $f$ over $X$. 

We will examine the problem of finding $\bm{x}^\star \in X^\star$ under the following conditions. 

\begin{enumerate}
\item[(C1)] There is an independent and identically distributed sample $\bm{\xi}_0, \bm{\xi}_1, \ldots$ of realizations of the random vector $\bm{\xi}$;
\item[(C2)] There is an oracle which, for a given input point $(\bm{x},\bm{\xi}) \in \mathbb{R}^d \times \Xi$, returns a stochastic gradient $\mathsf{G}(\bm{x},\bm{\xi})$ such that $\mathbb{E}[\mathsf{G}(\bm{x},\bm{\xi})] = \nabla f (\bm{x})$;
\item[(C3)] There exists a positive number $M$ such that, for all $\bm{x}\in X$, $\mathbb{E}[\|\mathsf{G}(\bm{x},\bm{\xi})\|^2] \leq M^2$.
\end{enumerate}

\section{Proposed Algorithm}
\label{sec:3}
This section describes the following algorithm (Algorithm \ref{algo:1}) for solving the problem of finding $\bm{x}^\star \in X^\star$ under (C1)--(C3).

\begin{algorithm} 
\caption{Adaptive learning rate optimization algorithm for solving the problem of finding $\bm{x}^\star \in X^\star$} 
\label{algo:1} 
\begin{algorithmic}[1] 
\REQUIRE
$(\alpha_n)_{n\in\mathbb{N}} \subset (0,1)$, $(\beta_n)_{n\in\mathbb{N}} \subset [0,1)$, 
$\gamma \in [0,1)$
\STATE
$n \gets 0$, $\bm{x}_0, \bm{m}_{-1} \in \mathbb{R}^d$, 
$\mathsf{H}_0 \in \mathbb{S}_{++}^d \cap \mathbb{D}^d$
\LOOP 
\STATE 
$\bm{m}_n := \beta_n \bm{m}_{n-1} + (1-\beta_n) \mathsf{G}(\bm{x}_n, \bm{\xi}_n)$
\STATE
$\displaystyle{\hat{\bm{m}}_n := \frac{\bm{m}_n}{1-\gamma^{n+1}}}$
\STATE
$\mathsf{H}_n \in \mathbb{S}_{++}^d \cap \mathbb{D}^d$
\STATE 
Find $\bm{\mathsf{d}}_n \in \mathbb{R}^d$ that solves $\mathsf{H}_n \bm{\mathsf{d}} = - \hat{\bm{m}}_n$
\STATE 
$\bm{x}_{n+1} := P_{X, \mathsf{H}_n} (\bm{x}_n + \alpha_n \bm{\mathsf{d}}_n)$
\STATE $n \gets n+1$
\ENDLOOP 
\end{algorithmic}
\end{algorithm}

Since $\mathsf{H}_n := \mathsf{diag}(h_{n,i})$ ($h_{n,i} > 0$) implies that there exists $\mathsf{H}_n^{-1} = \mathsf{diag}(h_{n,i}^{-1})$, step 7 in Algorithm \ref{algo:1} can be expressed as 
\begin{align}\label{7}
\bm{x}_{n+1} = P_{X, \mathsf{H}_n} 
\left[
\left( x_{n,i} - \frac{\alpha_n}{(1 - \gamma^{n+1}) h_{n,i}}
m_{n,i} \right)_{i=1}^d
\right].
\end{align}
We can see that Algorithm \ref{algo:1} adapts the learning rate $\alpha_n/((1 - \gamma^{n+1}) h_{n,i})$ for each $n\in\mathbb{N}$ and each $i=1,2,\ldots,d$. Throughout this paper, we call the parameters $\alpha_n$ and $\beta_n$ {\em sub-learning rates} for the learning rate $\alpha_n/((1 - \gamma^{n+1}) h_{n,i})$.

The convergence analyses of Algorithm \ref{algo:1} assume the following conditions.

\textbf{Assumption:} The sequence $(\mathsf{H}_n)_{n\in\mathbb{N}} \subset \mathbb{S}_{++}^d \cap \mathbb{D}^d$, denoted by $\mathsf{H}_n := \mathsf{diag}(h_{n,i})$, in Algorithm \ref{algo:1} satisfies the following conditions:
\begin{enumerate}
\item[(A3)] $h_{n+1,i} \geq h_{n,i}$ almost surely for all $n\in\mathbb{N}$ and all $i=1,2,\ldots,d$;
\item[(A4)] For all $i=1,2,\ldots,d$, a positive number $B_i$ exists such that $\sup \{ \mathbb{E}[h_{n,i}] \colon n \in \mathbb{N} \} \leq B_i$.
\end{enumerate}
Moreover,
\begin{enumerate}
\item[(A5)] $D := \max_{i=1,2,\ldots,d} \sup \{ (x_{i} - y_i)^2 \colon (x_i), (y_i) \in X \} < + \infty$.
\end{enumerate}
Assumption (A5) holds under the boundedness condition of $X$, which is assumed in \cite[p.1574]{nem2009} and \cite[p.2]{reddi2018}. Here, we provide some examples of $(\mathsf{H}_n)_{n\in\mathbb{N}}$ satisfying (A3) and (A4) when $X$ is bounded (i.e., (A5) holds). First, we consider $\mathsf{H}_n$ and $\bm{v}_n$ ($n\in\mathbb{N}$) defined for all $n\in\mathbb{N}$ by 
\begin{align}\label{adam}
\begin{split}
&\bm{v}_n := \delta \bm{v}_{n-1} + (1- \delta) \mathsf{G}(\bm{x}_n, \bm{\xi}_n) \odot \mathsf{G}(\bm{x}_n, \bm{\xi}_n),\\
&\bar{\bm{v}}_n := \frac{\bm{v}_n}{1 - \delta^{n+1}},\\
&\hat{\bm{v}}_n = (\hat{v}_{n,i}) := \left(\max \{ \hat{v}_{n-1,i}, \bar{v}_{n,i} \}\right),\\
&\mathsf{H}_n := \mathsf{diag} \left(\sqrt{\hat{v}_{n,i}} \right),
\end{split}
\end{align}
where $\bm{v}_{-1} = \hat{\bm{v}}_{-1} = \bm{0} \in \mathbb{R}^d$ and 
$\delta \in [0,1)$. Algorithm \ref{algo:1} with \eqref{adam} is based on the Adam algorithm\footnote{Adam uses $\mathsf{H}_n = \mathsf{diag}(\bar{v}_{n,i}^{1/2})$. We use $\hat{\bm{v}}_n = (\hat{v}_{n,i}) := (\max \{ \hat{v}_{n-1,i}, \bar{v}_{n,i} \})$ in \eqref{adam} so as to satisfy (A3). The modification of $\mathsf{H}_n$ defined by $\mathsf{diag} (\hat{v}_{n,i}^{1/2} + \epsilon)$ guarantees that $h_{n,i} \neq 0$, where $\epsilon > 0$ \cite{adam}.} \cite{adam}. The definitions of $\hat{\bm{v}}_n$ and $\mathsf{H}_n = \mathsf{diag}(h_{n,i}) = \mathsf{diag}(\hat{v}_{n,i}^{1/2}) \in \mathbb{S}_{++}^d \cap \mathbb{D}^d$ in \eqref{adam} obviously satisfy (A3). Step 7 in Algorithm \ref{algo:1} implies that $(\bm{x}_n)_{n\in\mathbb{N}} \subset X$. Accordingly, the boundedness of $X$ and (A2) ensure that $(\mathsf{G}(\bm{x}_n,\bm{\xi}_n))_{n\in\mathbb{N}}$ is almost surely bounded, i.e.,
\begin{align*}
M_1 := \sup \left\{ \left\| \mathsf{G}(\bm{x}_n,\bm{\xi}_n) \odot \mathsf{G}(\bm{x}_n,\bm{\xi}_n) \right\| \colon n\in\mathbb{N} \right\} < + \infty.
\end{align*}
Moreover, from the definition of $\bm{v}_n$ and the triangle inequality, we have, for all $n\in\mathbb{N}$,
\begin{align*}
\left\| \bm{v}_n \right\| 
&\leq \delta \left\| \bm{v}_{n-1} \right\| 
 + (1- \delta) M_1.
\end{align*}
Induction thus shows that, for all $n\in \mathbb{N}$, $\| \bm{v}_n \| = (\sum_{i=1}^d |v_{n,i}|^2)^{1/2} \leq M_1$, almost surely, which, together with the definition of $\bar{\bm{v}}_n$, implies that $\| \bar{\bm{v}}_n \| = (\sum_{i=1}^d |\bar{v}_{n,i}|^2)^{1/2} \leq M_1/(1-\delta)$. Accordingly, we have, for all $n\in\mathbb{N}$ and all $i=1,2,\ldots,d$,
\begin{align*}
|v_{n,i}|^2, |\bar{v}_{n,i}|^2 \leq \frac{M_1^2}{(1-\delta)^2}.
\end{align*}
The definition of $\hat{\bm{v}}_n$ and $\hat{\bm{v}}_{-1} = \bm{0}$ ensure that, for all $n\in\mathbb{N}$ and all $i=1,2,\ldots,d$, 
\begin{align*}
\mathbb{E}[h_{n,i}] := 
\mathbb{E}\left[\sqrt{\hat{v}_{n,i}} \right] \leq \frac{M_1}{1-\delta},
\end{align*}
which implies that (A4) holds.

Next, we consider $\mathsf{H}_n$ and $\bm{v}_n$ ($n\in\mathbb{N}$) defined for all $n\in\mathbb{N}$ by 
\begin{align}\label{amsgrad}
\begin{split}
&\bm{v}_n := \delta \bm{v}_{n-1} + (1- \delta) \mathsf{G}(\bm{x}_n, \bm{\xi}_n) \odot \mathsf{G}(\bm{x}_n, \bm{\xi}_n),\\
&\hat{\bm{v}}_n = (\hat{v}_{n,i}) := \left(\max \{ \hat{v}_{n-1,i}, v_{n,i} \}\right),\\
&\mathsf{H}_n := \mathsf{diag} \left(\sqrt{\hat{v}_{n,i}} \right),
\end{split}
\end{align}
where $\bm{v}_{-1} = \hat{\bm{v}}_{-1} = \bm{0} \in \mathbb{R}^d$ and 
$\delta \in [0,1)$. Algorithm \ref{algo:1} with \eqref{amsgrad} is the AMSGrad algorithm \cite{reddi2018}. A discussion similar to the one showing that $\mathsf{H}_n$ and $\bm{v}_n$ defined by \eqref{adam} satisfy (A3) and (A4) ensures that $\mathsf{H}_n$ and $\bm{v}_n$ defined by \eqref{amsgrad} satisfy (A3) and (A4); i.e., for all $n\in\mathbb{N}$ and all $i=1,2,\ldots,d$, 
\begin{align*}
\mathbb{E}[h_{n,i}] := \mathbb{E}\left[\sqrt{\hat{v}_{n,i}} \right] \leq M_1.
\end{align*}

\subsection{Constant sub-learning rate rule}
\label{subsec:3.1}
The following is the convergence analysis of Algorithm \ref{algo:1} with constant sub-learning rates. The proof of Theorem \ref{thm:1} is given in Appendix \ref{appen:1}.

\begin{thm}\label{thm:1}
Suppose that (A1)--(A5) and (C1)--(C3) hold and $(\bm{x}_n)_{n\in\mathbb{N}}$ is the sequence generated by Algorithm \ref{algo:1} with $\alpha_n := \alpha$ and $\beta_n := \beta$ ($n\in\mathbb{N}$). Then, for all $\bm{x} \in X$,
\begin{align*}
\limsup_{n \to + \infty} \mathbb{E}\left[ \langle \bm{x} - \bm{x}_n,
\nabla f (\bm{x}_n) \rangle \right] 
\geq - \frac{\tilde{B}^2 \tilde{M}^2}{2 \tilde{b}\tilde{\gamma}^2} \alpha 
- \frac{\tilde{M}\sqrt{Dd}}{\tilde{b}\tilde{\gamma}} \beta,
\end{align*}
where $\tilde{\gamma} := 1 - \gamma$, $\tilde{b} := 1 - \beta$, $\tilde{M}^2 := \max\{ \|\bm{m}_{-1}\|^2, M^2 \}$, $D$ is defined as in (A5), and 
$\tilde{B} := \sup\{ \max_{i=1,2,\ldots,d} h_{n,i}^{-1/2} \colon n\in\mathbb{N}\} < + \infty$.
\end{thm}

Algorithm \ref{algo:1} with $\alpha_n := \alpha$ and $\beta_n := \beta$ is as follows (see also \eqref{7}):
\begin{align}\label{7_1}
\begin{split}
&m_{n,i} = \beta m_{n-1,i} + (1-\beta) \mathsf{g}_{n,i},\\
&\bm{x}_{n+1} = P_{X, \mathsf{H}_n} 
\left[
\left( x_{n,i} - \frac{\alpha}{(1 - \gamma^{n+1}) h_{n,i}}
m_{n,i} \right)_{i=1}^d
\right],
\end{split}
\end{align}
where $(\mathsf{g}_{n,i}) := \mathsf{G}(\bm{x}_n, \bm{\xi}_n)$. Assumptions (A3) and (A4) guarantee that $(\mathbb{E}[h_{n,i}^{-1}])_{n\in\mathbb{N}}$ converges to $h_i^\star > 0$. Accordingly, the sequence of learning rates $({\mathbb{E}}[\alpha/((1 - {\gamma}^{n+1}) h_{n,i})])_{n\in\mathbb{N}}$ converges to $\alpha h_i^\star > 0$; i.e., the sequence of learning rates does not converge to zero. Therefore, we can see that Algorithm \ref{algo:1} with constant sub-learning rates is implementable in practice. Theorem \ref{thm:1} indicates that Algorithm \ref{algo:1} with small constant sub-learning rates, 
\begin{align}\label{Rate}
\alpha = \frac{1}{10^{a_1}}
\text{ and }
\beta = \frac{1}{10^{a_2}}
\quad (a_1,a_2 > 0),
\end{align}
approximates {$x^\star \in X^\star$} in the sense of {the existence of a positive real number $B$ such that 
\begin{align}\label{rate:1}
&\limsup_{n \to + \infty} \mathbb{E}\left[ \langle \bm{x} - \bm{x}_n,
\nabla f (\bm{x}_n) \rangle \right] 
\geq -B \left(\frac{1}{10^{a_1}} + \frac{1}{10^{a_2}} \right).
\end{align}
}
The previously reported results in \cite{adam} and \cite{reddi2018} used a fixed parameter $\beta := 0.9$ {and was studied under the convexity condition of $f$}. Meanwhile, {Theorem \ref{thm:1}}, \eqref{Rate}, {and \eqref{rate:1}} show that using a small constant sub-learning rate $\beta$ is an appropriate way {to find $\bm{x}^\star \in X^\star$}. { Theorem \ref{thm:1} leads to the following.
\begin{prop}\label{prop:0}
Suppose that (A1)--(A5) and (C1)--(C3) hold, $F(\cdot,\bm{\xi})$ is convex for almost every $\bm{\xi} \in \Xi$, and $(\bm{x}_n)_{n\in\mathbb{N}}$ is the sequence generated by Algorithm \ref{algo:1} with $\alpha_n := \alpha$ and $\beta_n := \beta$ ($n\in\mathbb{N}$). Then, 
\begin{align*}
\liminf_{n \to + \infty} \mathbb{E} \left[ f(\bm{x}_n) - f^\star \right] 
\leq
\frac{\tilde{B}^2 \tilde{M}^2}{2 \tilde{b}\tilde{\gamma}^2} \alpha 
+ \frac{\tilde{M}\sqrt{Dd}}{\tilde{b}\tilde{\gamma}} \beta,
\end{align*}
where $f^\star$ denotes the optimal value of the problem of minimizing $f$ over $X$, and $\tilde{\gamma}$, $\tilde{b}$, $\tilde{M}$, $D$, and $\tilde{B}$ are defined as in Theorem \ref{thm:1}.
\end{prop}
}

\subsection{Diminishing sub-learning rate rule}
\label{subsec:3.2}
The following is the convergence analysis of Algorithm \ref{algo:1} with diminishing sub-learning rates. The proof of Theorem \ref{thm:2} is given in Appendix \ref{appen:1}.

\begin{thm}\label{thm:2}
Suppose that (A1)--(A5) and (C1)--(C3) hold {and $(\bm{x}_n)_{n\in\mathbb{N}}$ is the sequence generated by Algorithm \ref{algo:1} with $\alpha_n$ and $\beta_n$ ($n\in\mathbb{N}$)}\footnote{ {The sub-learning rates $\alpha_n := 1/n^\eta$ ($\eta \in (1/2,1]$) and $\beta_n := \lambda^n$ ($\lambda \in (0,1)$) satisfy $\sum_{n=0}^{+\infty} \alpha_n = + \infty$, $\sum_{n=0}^{+\infty} \alpha_n^2 < + \infty$, and $\sum_{n=0}^{+\infty} \alpha_n \beta_n < + \infty$ (by $\lim_{n\to + \infty} (\alpha_{n+1}\beta_{n+1})/(\alpha_n \beta_n) = \lambda \in (0,1)$).}} {satisfying $\sum_{n=0}^{+\infty} \alpha_n = + \infty$, $\sum_{n=0}^{+\infty} \alpha_n^2 < + \infty$, and $\sum_{n=0}^{+\infty} \alpha_n \beta_n < + \infty$. Then, for all $\bm{x} \in X$,
\begin{align}\label{sup}
\limsup_{n \to + \infty} \mathbb{E}\left[ \langle \bm{x} - \bm{x}_n,
\nabla f (\bm{x}_n) \rangle \right] 
\geq 0.
\end{align}
Moreover, if $\alpha_n := 1/n^\eta$ ($\eta \in [1/2,1)$)}\footnote{{Algorithm \ref{algo:1} with $\alpha_n := 1/\sqrt{n}$ and $\beta_n := \lambda^n$ does not satisfy \eqref{sup}. However, Algorithm \ref{algo:1} with $\alpha_n := 1/\sqrt{n}$ and $\beta_n := \lambda^n$ achieves a convergence rate of $\mathcal{O}(1/\sqrt{n})$.}} {and if $\beta_n := \lambda^n$ ($\lambda \in (0,1)$), then Algorithm \ref{algo:1} achieves the following convergence rate:
\begin{align*}
\frac{1}{n} \sum_{k=1}^n \mathbb{E}\left[ \langle \bm{x} - \bm{x}_k,
\nabla f (\bm{x}_k) \rangle \right]
\geq 
- \mathcal{O} \left( \frac{1}{n^{1 - \eta}} \right).
\end{align*}
}
\end{thm}

Algorithm \ref{algo:1} with diminishing sub-learning rates $\alpha_n$ and $\beta_n$ is as follows (see also \eqref{7}):
\begin{align}\label{7_2}
\begin{split}
&m_{n,i} = \beta_n m_{n-1,i} + (1-\beta_n) \mathsf{g}_{n,i},\\
&\bm{x}_{n+1} = P_{X, \mathsf{H}_n} 
\left[
\left( x_{n,i} - \frac{\alpha_n}{(1 - {\gamma}^{n+1}) h_{n,i}}
m_{n,i} \right)_{i=1}^d
\right].
\end{split}
\end{align}
Assumptions (A3) and (A4) {and $\lim_{n\to + \infty} \alpha_n = 0$} guarantee that the sequence of learning rates $({\mathbb{E}}[\alpha_n/((1 - {\gamma}^{n+1}) h_{n,i})])_{n\in\mathbb{N}}$ converges to zero. This means that Algorithm \ref{algo:1} with diminishing sub-learning rates would not be implementable in practice. However, {\eqref{sup} in} Theorem \ref{thm:2} guarantees the convergence of Algorithm \ref{algo:1} with diminishing sub-learning rates to {a point in $X^\star$ in the sense of the existence of an accumulation point of $(\bm{x}_n)_{n\in\mathbb{N}}$ belonging to $X^\star$}. 
{The following proposition enables us to compare Algorithm \ref{algo:1} with Adam \cite{adam} and AMSGrad \cite{reddi2018}.}
{
\begin{prop}\label{prop:1}
Suppose that (A1)--(A5) and (C1)--(C3) hold, $F(\cdot,\bm{\xi})$ is convex for almost every $\bm{\xi} \in \Xi$, and $(\bm{x}_n)_{n\in\mathbb{N}}$ is the sequence generated by Algorithm \ref{algo:1} with $\alpha_n := 1/n^\eta$ and $\beta_n := \lambda^n$ ($n\in\mathbb{N}$), where $\eta \in [1/2,1]$ and $\lambda \in (0,1)$. Then, under $\eta \in (1/2,1]$, 
\begin{align*}
\liminf_{n \to + \infty} \mathbb{E} \left[ f(\bm{x}_n) - f^\star \right] = 0,
\end{align*}
where $f^\star$ denotes the optimal value of the problem of minimizing $f$ over $X$. Moreover, under $\eta \in [1/2,1)$, any accumulation point of $(\tilde{\bm{x}}_n)_{n\in\mathbb{N}}$ defined by $\tilde{\bm{x}}_n := (1/n) \sum_{k=1}^n \bm{x}_k$ almost surely belongs to $X^\star$, and Algorithm \ref{algo:1} achieves the following convergence rate:
\begin{align*}
\mathbb{E}\left[ f(\tilde{\bm{x}}_n) - f^\star \right]
= 
\mathcal{O} \left( \frac{1}{n^{1 - \eta}} \right).
\end{align*}
\end{prop}
}

\subsection{Comparison of Algorithm \ref{algo:1} with the Existing Algorithms}
\label{sec:4}
The main objective of adaptive learning rate optimization algorithms is to solve Problem with $f(\bm{x}) = \mathbb{E}[F(\bm{x}, \xi)] = (1/T) \sum_{t=1}^T f_t(\bm{x})$ under (A1)--(A2) and (C1)--(C3), i.e.,
\begin{align}\label{prob:2}
\text{minimize } \sum_{t=1}^T f_t (\bm{x}) \text{ subject to } \bm{x}\in X,
\end{align}
where $T$ is the number of training examples, $f_t (\cdot) = F(\cdot, t) \colon \mathbb{R}^d \to \mathbb{R}$ $(t=1,2,\ldots,T)$ is a differentiable, convex loss function, and $X \subset \mathbb{R}^d$ is bounded, closed, and convex (i.e., (A5) holds). The performance measure of adaptive learning rate optimization algorithms to solve problem \eqref{prob:2} is called the {\em regret} on a sequence of $(f_t(\bm{x}_t))_{t=1}^T$ defined as follows:
\begin{align}\label{REG}
R(T) := \sum_{t=1}^T f_t (\bm{x}_t) - f^\star,
\end{align} 
where {$f^\star$ denotes the optimal value of problem \eqref{prob:2} and} $(\bm{x}_t)_{t=1}^T \subset X$ is the sequence generated by a learning algorithm. Adam \cite{adam} is useful for training deep neural networks. Theorem 4.1 in \cite{adam} indicates that Adam ensures that there is a positive real number $D$ such that $R(T)/T \leq D/\sqrt{T}$. However, Theorem 1 in \cite{reddi2018} shows that a counter-example to Theorem 4.1 in \cite{adam} exists.

The difference between Adam and Algorithm \ref{algo:1} with \eqref{adam} is in the definitions of $\bm{m}_n$ and $\mathsf{H}_n$, {i.e., Adam is defined by using $\mathsf{h}_{n,i} = \sqrt{\bar{v}_{n,i}}$, while Algorithm \ref{algo:1} with \eqref{adam} uses $\mathsf{h}_{n,i} = \sqrt{\hat{v}_{n,i}}$.} {While Adam does not converge to a solution of problem \eqref{prob:2} \cite[Theorem 1]{reddi2018}, this difference leads to Proposition \ref{prop:0}}, {indicating that Algorithm \ref{algo:1} with \eqref{adam}}, $\alpha_n := \alpha$, and $\beta_n := \beta$ ensures {that there exist $C_1, C_2 > 0$ such that $\liminf_{n\to + \infty} \mathbb{E}[f(\bm{x}_n) - f^\star ] \leq C_1 \alpha + C_2 \beta$.} This implies that, if we can use sufficiently small constant learning rates $\alpha$ and $\beta$ (see, e.g., \eqref{Rate}), then {Algorithm \ref{algo:1} with \eqref{adam}} approximates the solution of problem \eqref{prob:2}. Although the previously reported results in \cite{adam,reddi2018} considered only the case where $(\alpha_t)_{t=1}^T$ is diminishing (e.g., $\alpha_t := \alpha/\sqrt{t}$), the above result from {Proposition \ref{prop:0}} guarantees that {Algorithm \ref{algo:1}} with {\em constant} sub-learning rates can solve problem \eqref{prob:2}. Moreover, {Proposition \ref{prop:1}} indicates that {Algorithm \ref{algo:1} with \eqref{adam}, $\alpha_n := 1/n^\eta$ ($\eta \in (1/2,1]$), and $\beta_n := \lambda^{n}$ ($\lambda \in (0,1)$) satisfies $\liminf_{n\to + \infty} \mathbb{E}[f(\bm{x}_n) - f^\star] = 0$ and} that any accumulation point of $(\tilde{\bm{x}}_n)_{n\in\mathbb{N}}$ generated by {Algorithm \ref{algo:1} with \eqref{adam}}, $\alpha_n := 1/n^\eta$ ($\eta \in [1/2,1)$), and $\beta_n := \lambda^{n}$ ($\lambda \in (0,1)$) almost surely belongs to the solution set of problem \eqref{prob:2}. This implies that {Algorithm \ref{algo:1}} with {\em diminishing} sub-learning rates can solve problem \eqref{prob:2}. 

AMSGrad \cite{reddi2018} was proposed as a way to guarantee convergence and preserve the practical benefits of Adam. The AMSGrad algorithm has the following property {(see Corollary 4.2 in \cite{adam} and Theorem 4 and Corollary 1 in \cite{reddi2018})}: Suppose that $\beta_{t} := {\nu} \lambda^{t}$ $({\nu}, \lambda \in (0,1))$, ${\theta} := {\nu}/\sqrt{{\delta}} < 1$, and $\alpha_t := \alpha/\sqrt{t}$ $(\alpha > 0)$. AMSGrad ensures that there is a positive real number $\hat{D}$ such that
\begin{align}\label{regret}
\frac{R(T)}{T} = 
\frac{1}{T} \left(\sum_{t=1}^T f_t (\bm{x}_t) - f^\star \right)
\leq \hat{D} \sqrt{\frac{1+\ln T}{T}}.
\end{align}
From the discussion on \eqref{7_2}, the sequences of the learning rates $(\alpha_t/((1-{\gamma}^{t+1})\sqrt{\bar{v}_{t,i}}))$ in Adam {(see also \eqref{adam})} and $(\alpha_t/\sqrt{\hat{v}_{t,i}})$ in AMSGrad {(see also \eqref{amsgrad})} converge to zero when $t$ diverges. Hence, Adam and AMSGrad would not be implementable in practice. 

{Algorithm \ref{algo:1} with \eqref{amsgrad} when} $n=1,2,\ldots,T$ coincides with AMSGrad. Although the previously reported results in \cite{reddi2018} (e.g., \eqref{regret}) show that the regret is low for a sufficiently large parameter $T$, but $(\bm{x}_t)_{t=1}^T$ defined by AMSGrad does not always approximate a solution of problem \eqref{prob:2}, as can be seen in \eqref{regret}. Meanwhile, {Proposition \ref{prop:0} indicates that Algorithm \ref{algo:1} with \eqref{amsgrad}, $\alpha_n := \alpha$, and $\beta_n := \beta$ approximates a solution of problem \eqref{prob:2} in the sense of $\liminf_{n\to + \infty} \mathbb{E}[f(\bm{x}_n) - f^\star ] \leq C_1 \alpha + C_2 \beta$}. {Proposition \ref{prop:1} indicates that Algorithm \ref{algo:1} with \eqref{adam}, $\alpha_n := 1/n^\eta$ ($\eta \in (1/2,1]$), and $\beta_n := \lambda^{n}$ ($\lambda \in (0,1)$) satisfies $\liminf_{n\to + \infty} \mathbb{E}[f(\bm{x}_n) - f^\star] = 0$ and that any accumulation point of $(\tilde{\bm{x}}_n)_{n\in\mathbb{N}}$ generated by Algorithm \ref{algo:1} with \eqref{amsgrad}, $\alpha_n := 1/n^\eta$ ($\eta \in [1/2,1)$), and $\beta_n := \lambda^{n}$ ($\lambda \in (0,1)$) almost surely belongs to the solution set of problem \eqref{prob:2}}. 

\section{{Numerical Experiments}}\label{sec:5}
We examined the behavior of Algorithm \ref{algo:1} with different sub-learning rates. The adaptive learning rate optimization algorithms with $\delta = 0.999$ \cite{adam,reddi2018} used in the experiments were as follows, where the initial point $\bm{x}_0$ was randomly chosen.

Algorithm \ref{algo:1} with constant sub-learning rates:
\begin{itemize}
\item ADAM-C1: Algorithm \ref{algo:1} with \eqref{adam}, $\gamma = 0.9$, $\alpha_n = 10^{-3}$, and $\beta_n = 0.9$
\item ADAM-C2: Algorithm \ref{algo:1} with \eqref{adam}, $\gamma = 0.9$, $\alpha_n = 10^{-3}$, and $\beta_n = 10^{-3}$
\item ADAM-C3: Algorithm \ref{algo:1} with \eqref{adam}, $\gamma = 0.9$, $\alpha_n = 10^{-2}$, and $\beta_n = 10^{-2}$
\item AMSG-C1: Algorithm \ref{algo:1} with \eqref{amsgrad}, $\gamma = 0$, $\alpha_n = 10^{-3}$, and $\beta_n = 0.9$
\item AMSG-C2: Algorithm \ref{algo:1} with \eqref{amsgrad}, $\gamma = 0$, $\alpha_n = 10^{-3}$, and $\beta_n = 10^{-3}$
\item AMSG-C3: Algorithm \ref{algo:1} with \eqref{amsgrad}, $\gamma = 0$, $\alpha_n = 10^{-2}$, and $\beta_n = 10^{-2}$
\item MAMSG-C1: Algorithm \ref{algo:1} with \eqref{amsgrad}, $\gamma = 0.1$, $\alpha_n = 10^{-3}$, and $\beta_n = 0.9$
\item MAMSG-C2: Algorithm \ref{algo:1} with \eqref{amsgrad}, $\gamma = 0.1$, $\alpha_n = 10^{-3}$, and $\beta_n = 10^{-3}$
\item MAMSG-C3: Algorithm \ref{algo:1} with \eqref{amsgrad}, $\gamma = 0.1$, $\alpha_n = 10^{-2}$, and $\beta_n = 10^{-2}$
\end{itemize}

Algorithm \ref{algo:1} with diminishing sub-learning rates:
\begin{itemize}
\item ADAM-D1: Algorithm \ref{algo:1} with \eqref{adam}, $\gamma = 0.9$, $\alpha_n = 1/\sqrt{n}$, and $\beta_n = 1/2^n$
\item ADAM-D2: Algorithm \ref{algo:1} with \eqref{adam}, $\gamma = 0.9$, $\alpha_n = 1/n^{3/4}$, and $\beta_n = 1/2^n$
\item ADAM-D3: Algorithm \ref{algo:1} with \eqref{adam}, $\gamma = 0.9$, $\alpha_n = 1/n$, and $\beta_n = 1/2^n$
\item AMSG-D1: Algorithm \ref{algo:1} with \eqref{amsgrad}, $\gamma = 0$, $\alpha_n = 1/\sqrt{n}$, and $\beta_n = 1/2^n$
\item AMSG-D2: Algorithm \ref{algo:1} with \eqref{amsgrad}, $\gamma = 0$, $\alpha_n = 1/n^{3/4}$, and $\beta_n = 1/2^n$
\item AMSG-D3: Algorithm \ref{algo:1} with \eqref{amsgrad}, $\gamma = 0$, $\alpha_n = 1/n$, and $\beta_n = 1/2^n$
\item MAMSG-D1: Algorithm \ref{algo:1} with \eqref{amsgrad}, $\gamma = 0.1$, $\alpha_n = 1/\sqrt{n}$, and $\beta_n = 1/2^n$
\item MAMSG-D2: Algorithm \ref{algo:1} with \eqref{amsgrad}, $\gamma = 0.1$, $\alpha_n = 1/n^{3/4}$, and $\beta_n = 1/2^n$
\item MAMSG-D3: Algorithm \ref{algo:1} with \eqref{amsgrad}, $\gamma = 0.1$, $\alpha_n = 1/n$, and $\beta_n = 1/2^n$
\end{itemize}

Although the parameters in ADAM-C1 are the same as those used in Adam \cite{adam}, ADAM-C1 (Algorithm \ref{algo:1} with \eqref{adam}) is a modification of Adam \cite{adam} that guarantees convergence. AMSG-C1 coincides with AMSGrad \cite{reddi2018}. We implemented ADAM-C$i$ (resp. AMSG-C$i$) ($i=2,3$) so that we could compare ADAM-C1 (resp. AMSG-C1) with the proposed algorithms with small constant sub-learning rates. We implemented ADAM-D$i$ and AMSG-D$i$ ($i = 1,2,3$), with diminishing sub-learning rates satisfying the conditions of Theorem \ref{thm:2}, and compared their performance with those of the algorithms with constant sub-learning rates. MAMSG-C$i$ (resp. MAMSG-D$i$) ($i=1,2,3$) with $\gamma = 0.1$ is a modification of AMSG-C$i$ (resp. AMSG-D$i$) with $\gamma = 0$.

The experiments used a fast scalar computation server at Meiji University. The environment has two Intel(R) Xeon(R) Gold 6148 (2.4 GHz, 20 cores) CPUs, an NVIDIA Tesla V100 (16GB, 900Gbps) GPU and a Red Hat Enterprise Linux 7.6 operating system. The experimental code was written in Python 3.8.2, and we used the NumPy 1.17.3 package and PyTorch 1.3.0 package.

\subsection{Text classification}
First, we considered a long short-term memory (LSTM) for text classification. The LSTM is an artificial recurrent neural network (RNN) architecture used in the field of deep learning for natural language processing. We used the IMDb dataset\footnote{\url{https://datasets.imdbws.com/}} for text classification tasks. The dataset contains 50,000 movie reviews along with their associated binary sentiment polarity labels. It is split  into 25,000 training and 25,000 test sets. A multilayer neural network was used to classify the dataset. The LSTM used in the experiment included  an affine layer and a sigmoid function as an activation function for the output. A loss function was the binary cross entropy (BCE). 

Figures \ref{fig:1_c} and \ref{fig:1_d} indicate that Algorithm \ref{algo:1} with constant sub-learning rates performed better than Algorithm \ref{algo:1} with diminishing ones in terms of both the train loss and accuracy score.

\subsection{Image classification}
Next, we considered a Residual Network (ResNet), which is a relatively deep model based on a convolutional neural network (CNN), for image classification. We used the CIFAR-10 dataset\footnote{\url{https://www.cs.toronto.edu/~kriz/cifar.html}}, which is a benchmark for image classification. The dataset consists of 60,000 color images ($32 \times 32$) in 10 classes, with 6,000 images per class. There are 50,000 training images and 10,000 test images. The test batch contained exactly 1,000 randomly selected images from each class. A 20-layer ResNet (ResNet-20) was organized into  19 convolutional layers that had $3 \times 3$ filters and a 10-way fully-connected layer with a softmax function. We used the cross entropy as the loss function for fitting ResNet in accordance with the common strategy in image classification.

Figures \ref{fig:2_c} and \ref{fig:2_d} indicate that Algorithm \ref{algo:1} with constant sub-learning rates performed better than Algorithm \ref{algo:1} with diminishing ones in terms of both the train loss and accuracy score. The reason why Algorithm \ref{algo:1} with diminishing sub-learning rates did not work so well is that the learning rates are approximately zero for a number of iterations. Meanwhile, Algorithm \ref{algo:1} with constant sub-learning rates worked well since its learning rates are not zero all of the time.

\begin{figure*}[htbp]
\vspace*{-1.5cm}
\hspace*{-3cm}
\subfloat[]{
\includegraphics[width=9truecm]{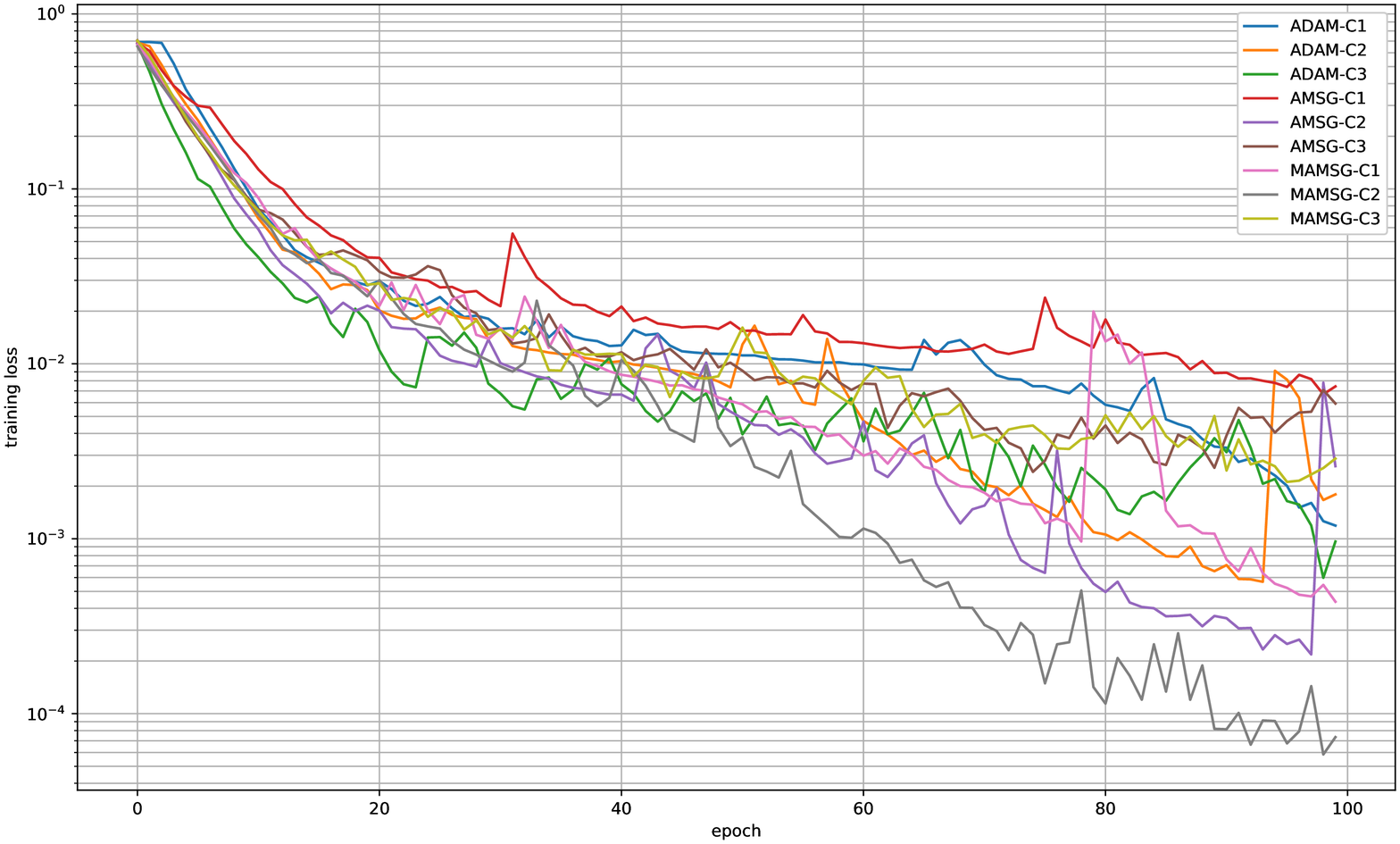}
\label{fig:1_c_l}
}
\hfil
\subfloat[]{
\includegraphics[width=9truecm]{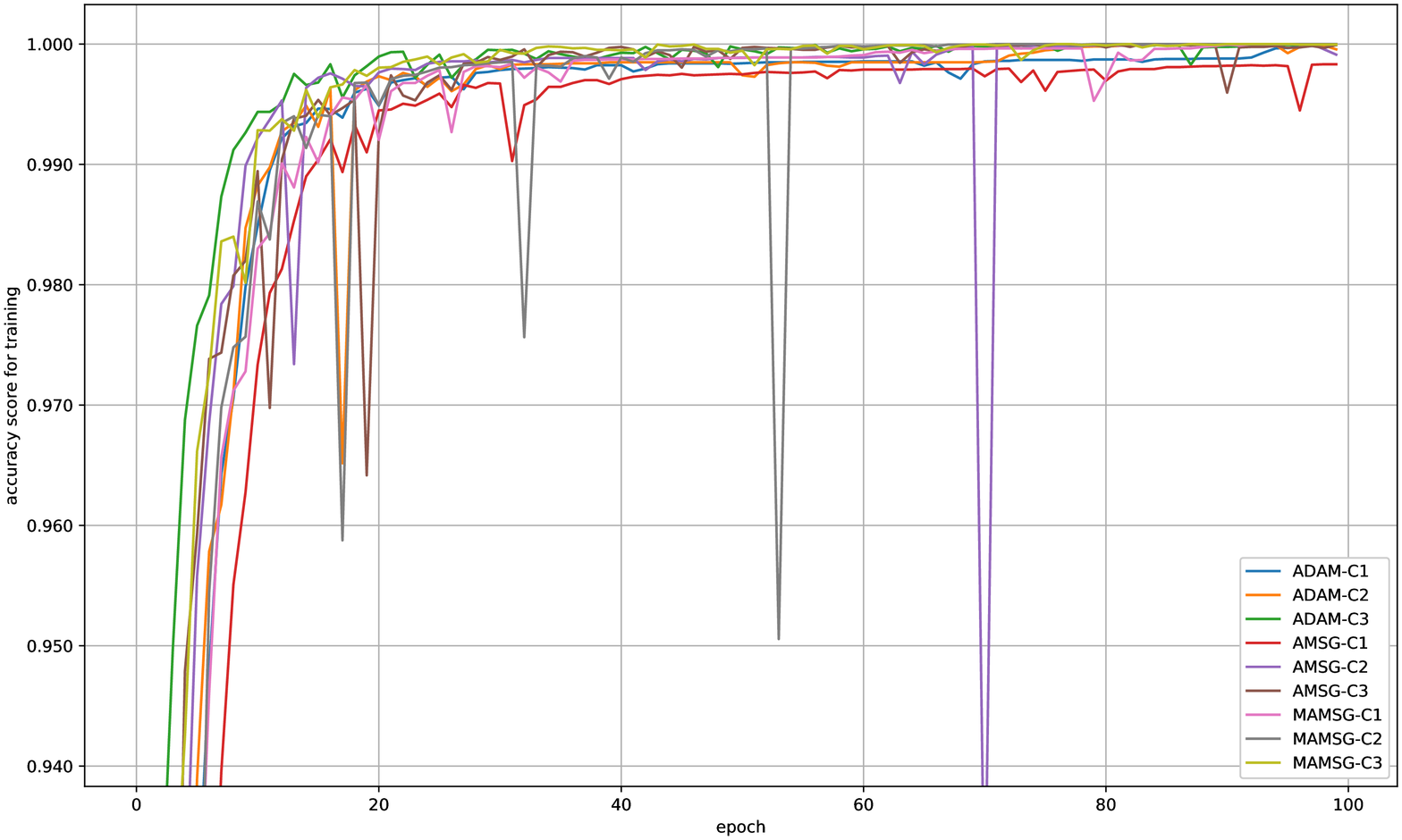}
\label{fig:1_c_a}
}
\hfil
\subfloat[]{
\includegraphics[width=9truecm]{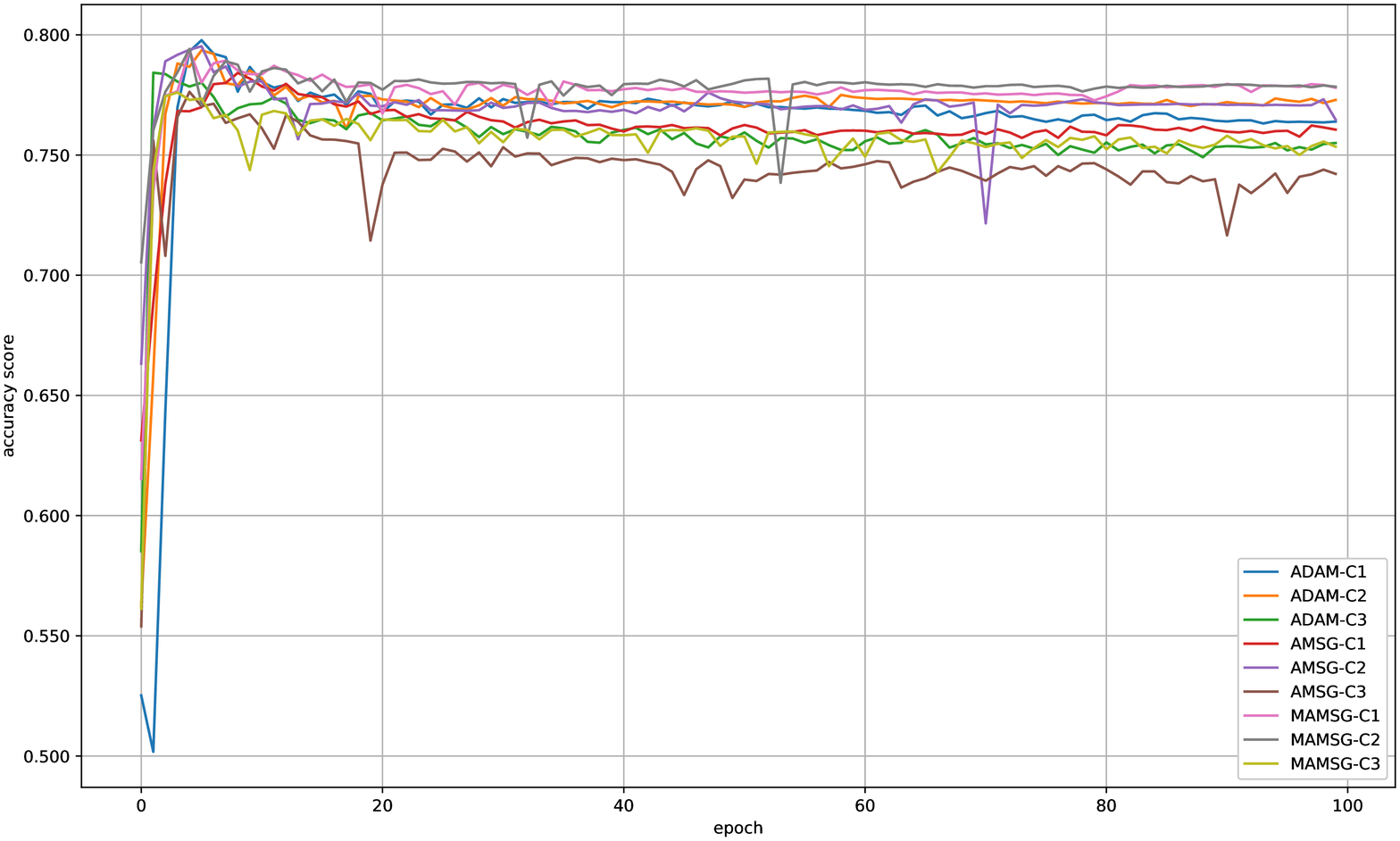}
\label{fig:1_c_a_t}
}

\caption{(a) Training loss function value, (b) training classification accuracy score, and (c) test classification accuracy score for Algorithm \ref{algo:1} with constant sub-learning rates versus number of epochs on the IMDb dataset}
\label{fig:1_c}
\end{figure*}

\begin{figure*}[htbp]
\vspace*{-1.5cm}
\hspace*{-3cm}
\subfloat[]{
\includegraphics[width=9truecm]{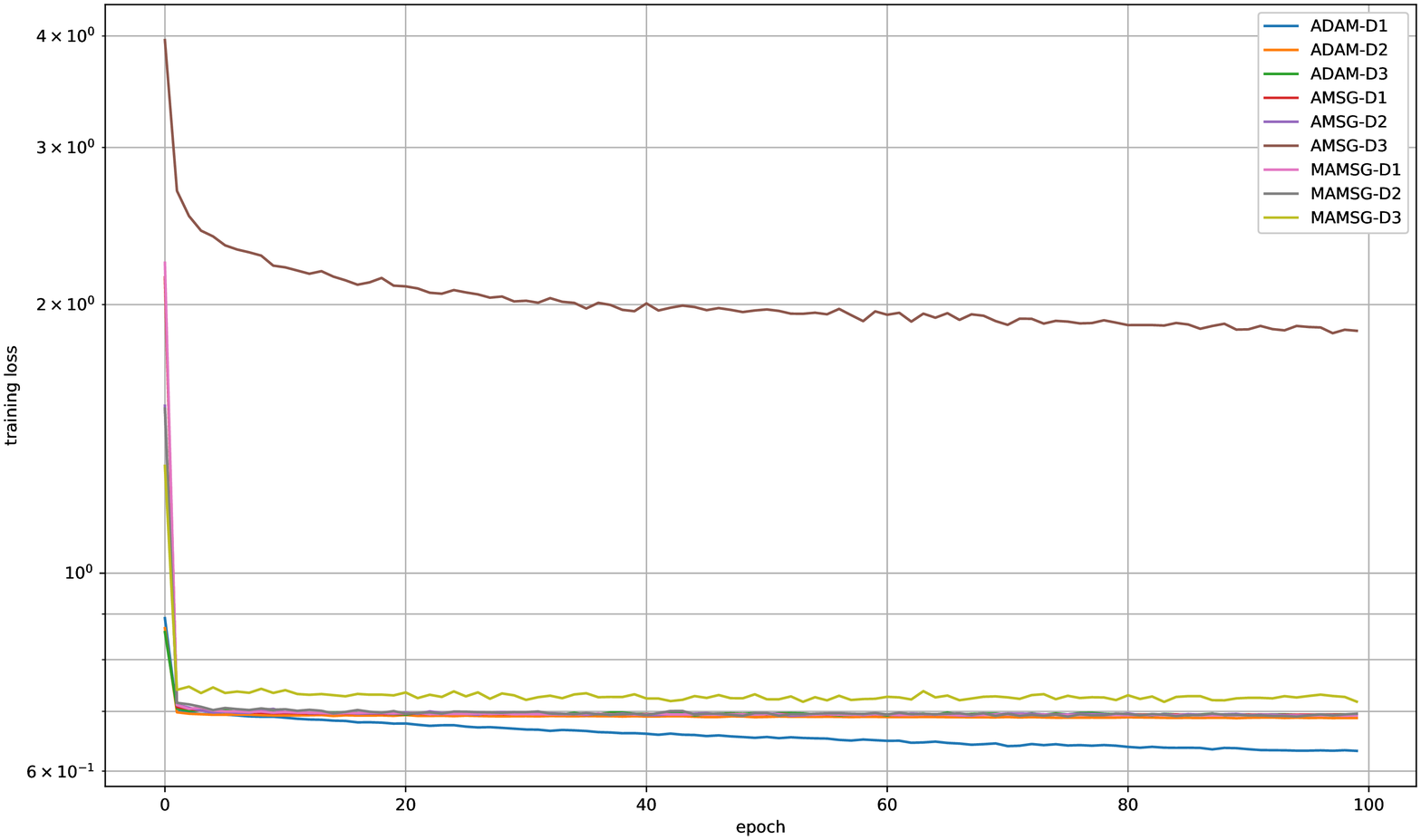}
\label{fig:1_d_l}
}
\hfil
\subfloat[]{
\includegraphics[width=9truecm]{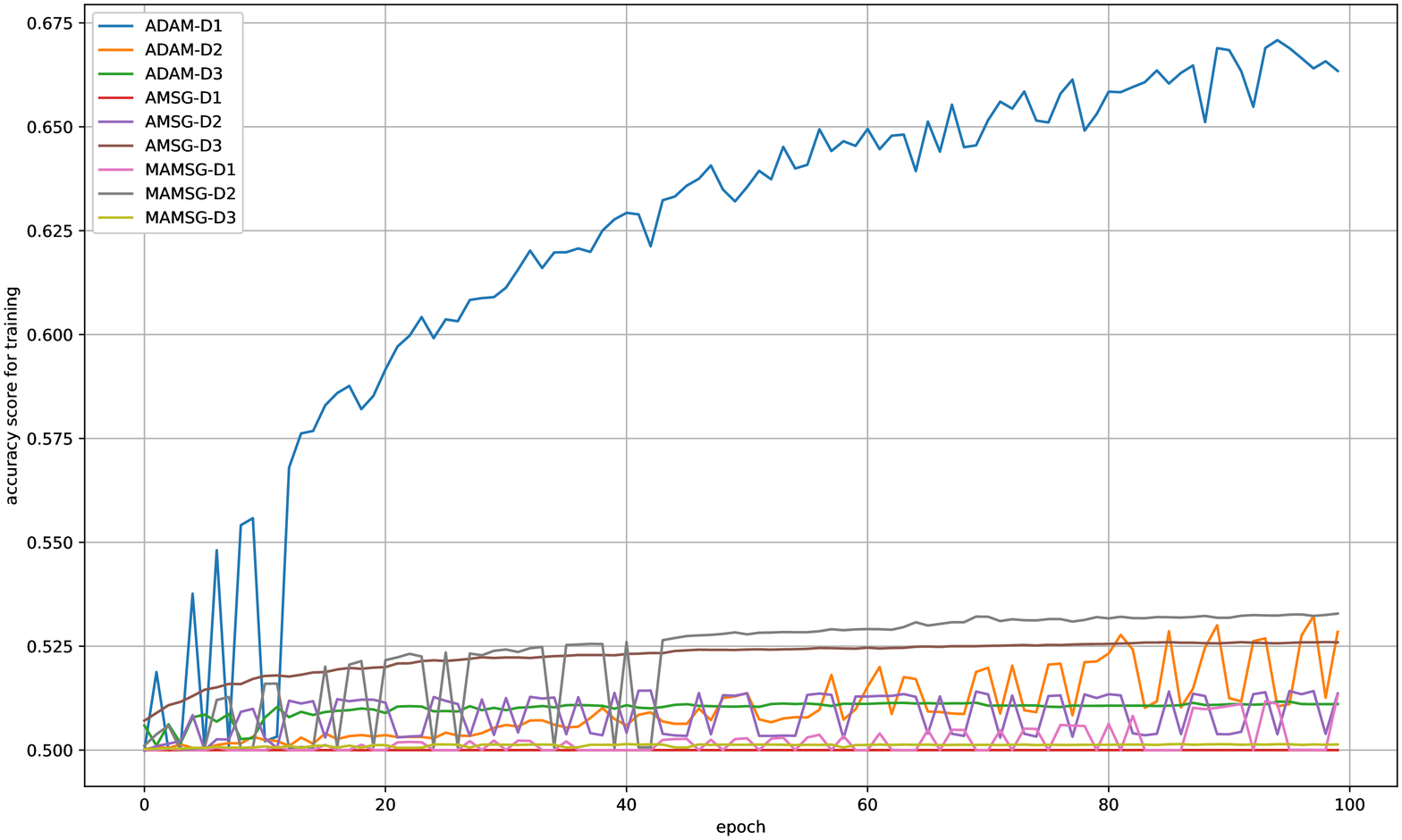}
\label{fig:1_d_a}
}
\hfil
\subfloat[]{
\includegraphics[width=9truecm]{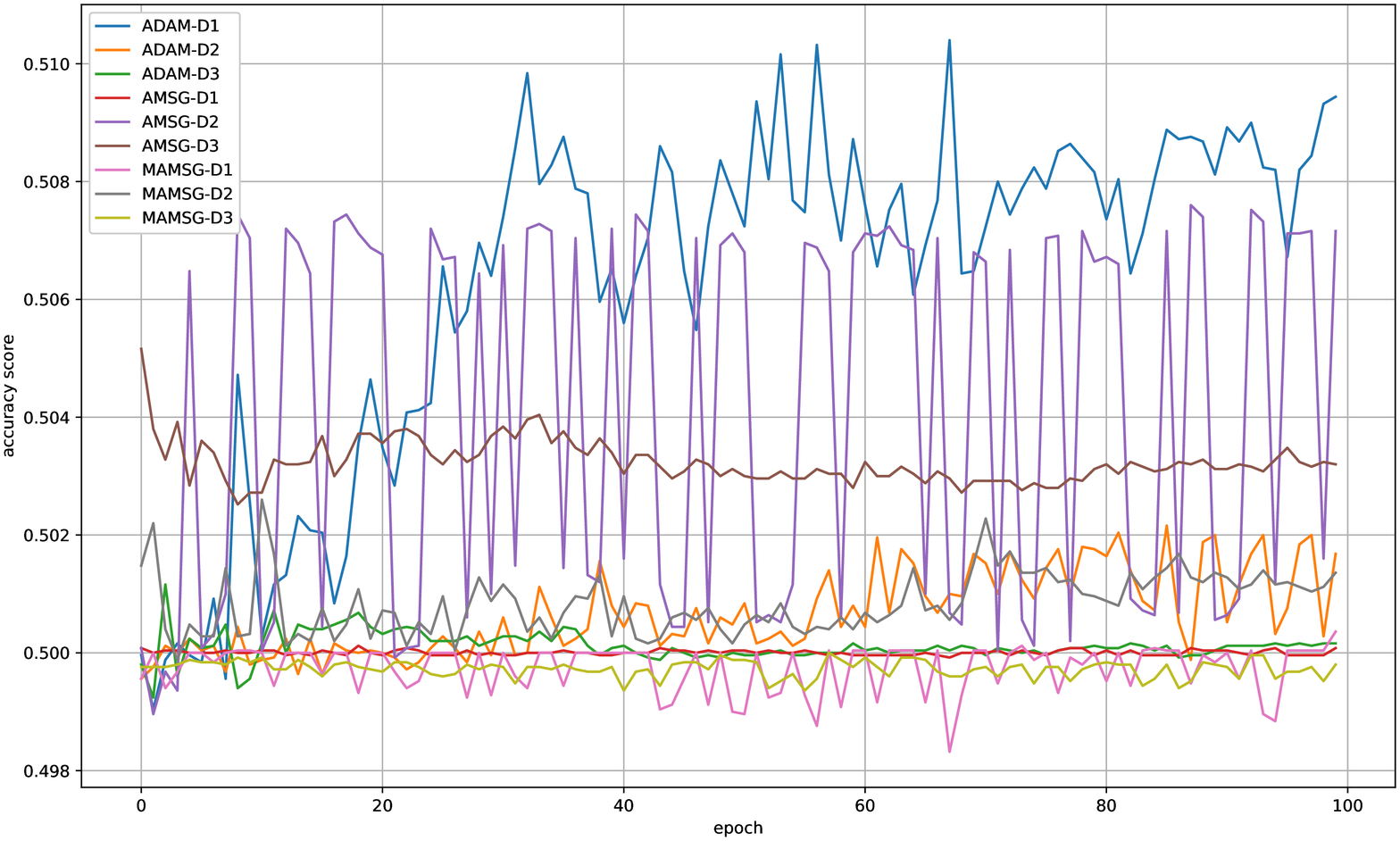}
\label{fig:1_d_a_t}
}
\caption{(a) Training loss function value, (b) training classification accuracy score, and (c) test classification accuracy score for Algorithm \ref{algo:1} with diminishing sub-learning rates versus number of epochs on the IMDb dataset}
\label{fig:1_d}
\end{figure*}

\begin{figure*}[htbp]
\vspace*{-1.5cm}
\hspace*{-3cm}
\subfloat[]{
\includegraphics[width=9truecm]{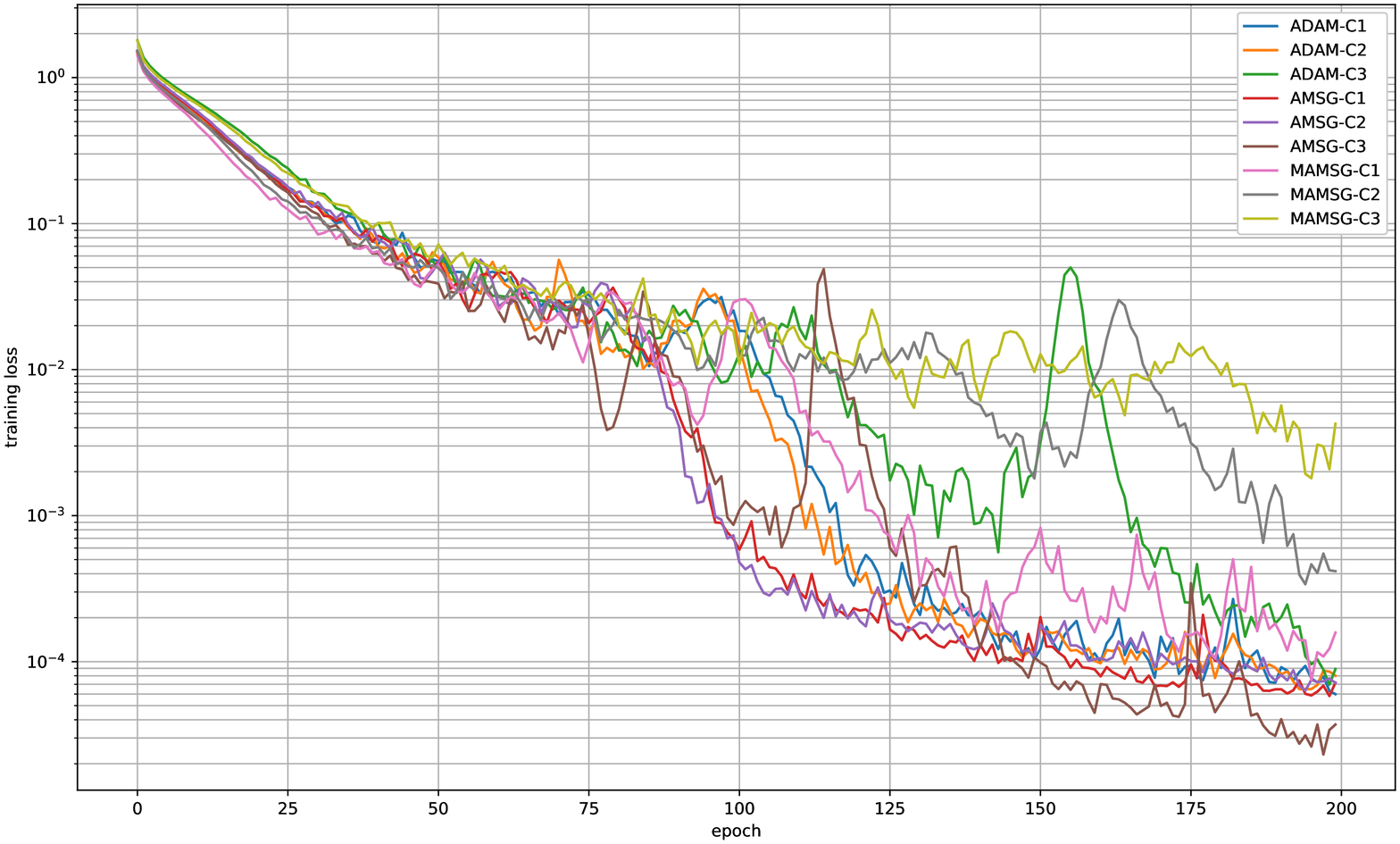}
\label{fig:2_c_l}
}
\hfil
\subfloat[]{
\includegraphics[width=9truecm]{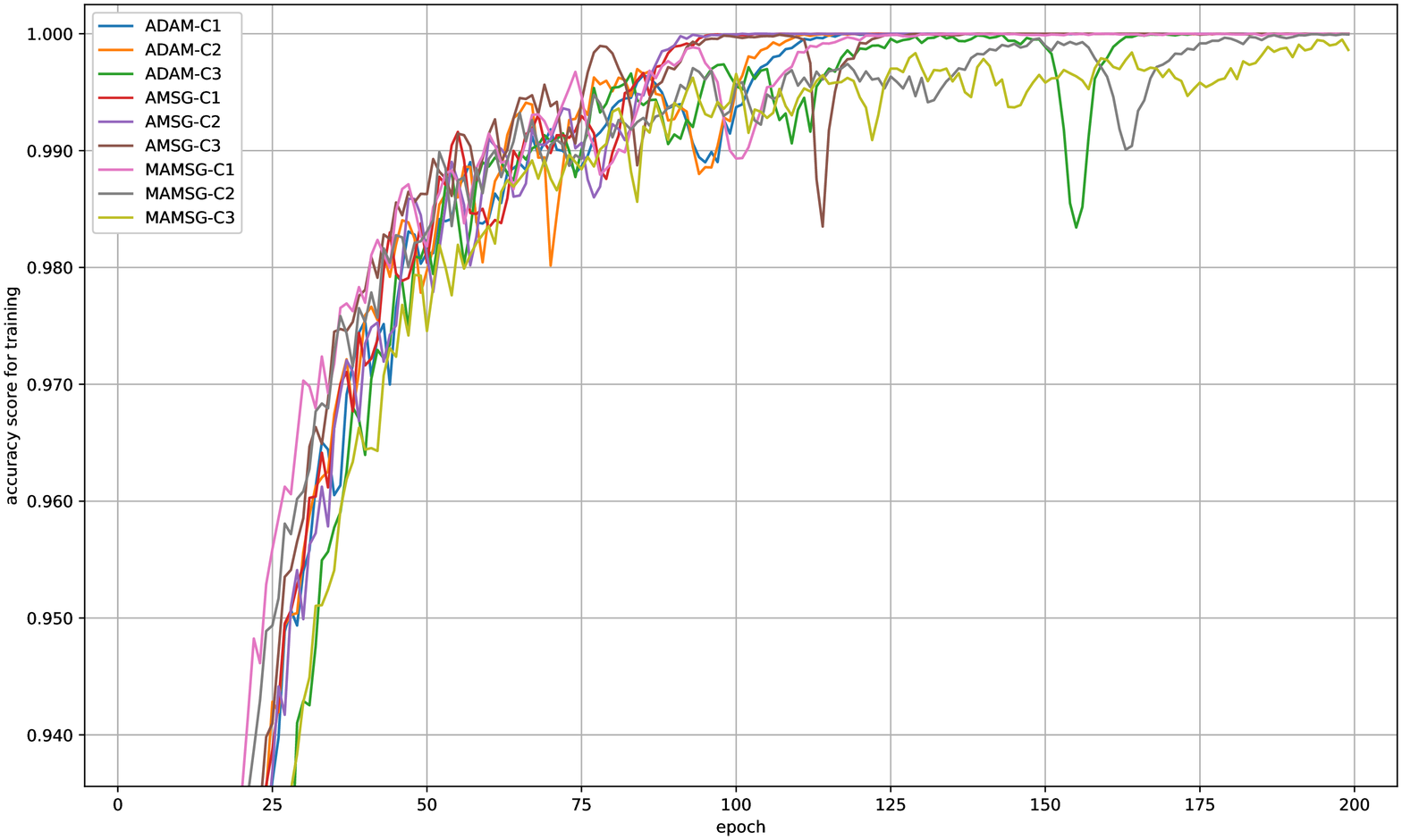}
\label{fig:2_c_a}
}
\hfil
\subfloat[]{
\includegraphics[width=9truecm]{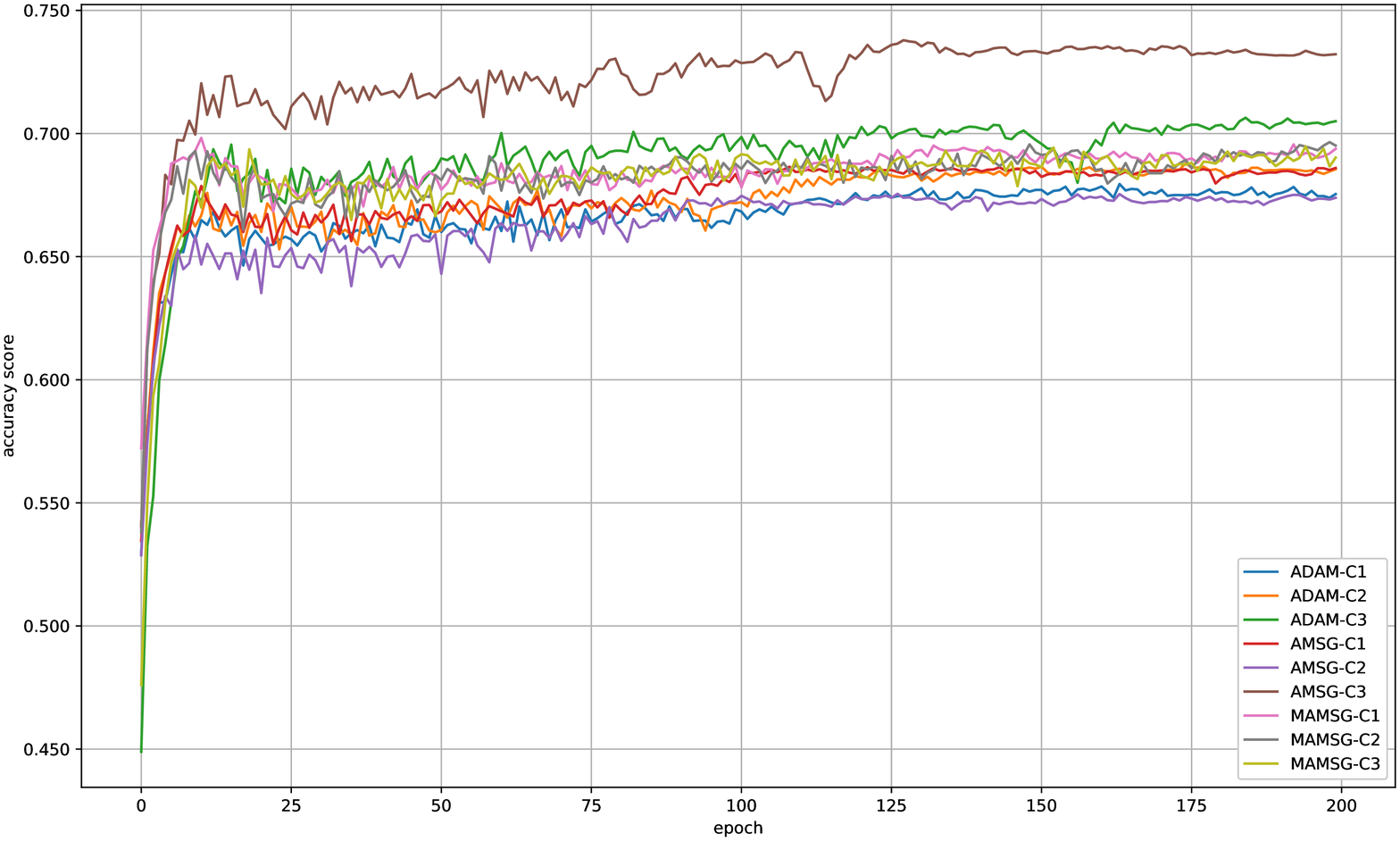}
\label{fig:2_c_a_t}
}
\caption{(a) Training loss function value, (b) training classification accuracy score, and (c) test classification accuracy score for Algorithm \ref{algo:1} with constant sub-learning rates versus number of epochs on the CIFAR-10 dataset}
\label{fig:2_c}
\end{figure*}

\begin{figure*}[htbp]
\vspace*{-1.5cm}
\hspace*{-3cm}
\subfloat[]{
\includegraphics[width=9truecm]{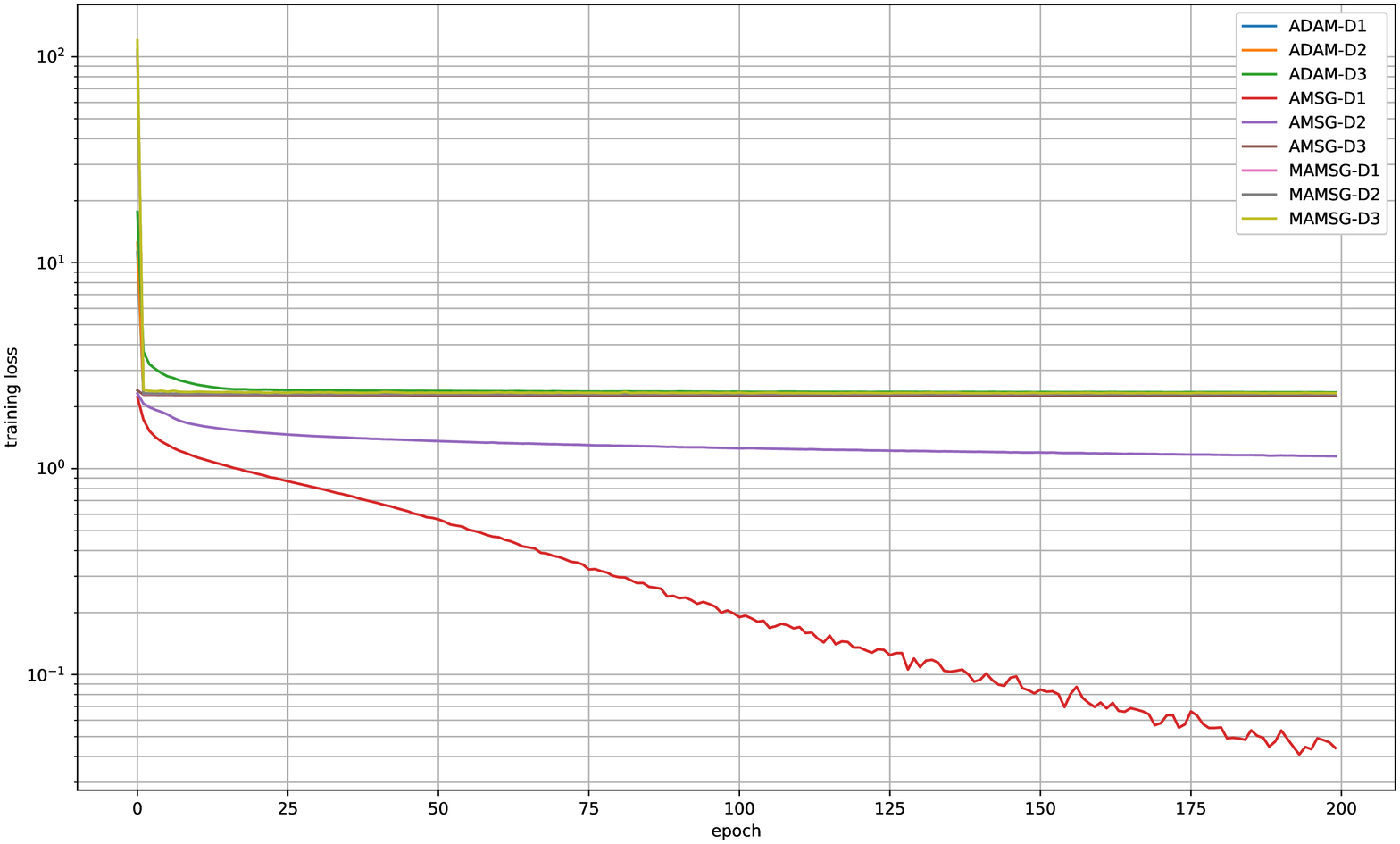}
\label{fig:2_d_l}
}
\hfil
\subfloat[]{
\includegraphics[width=9truecm]{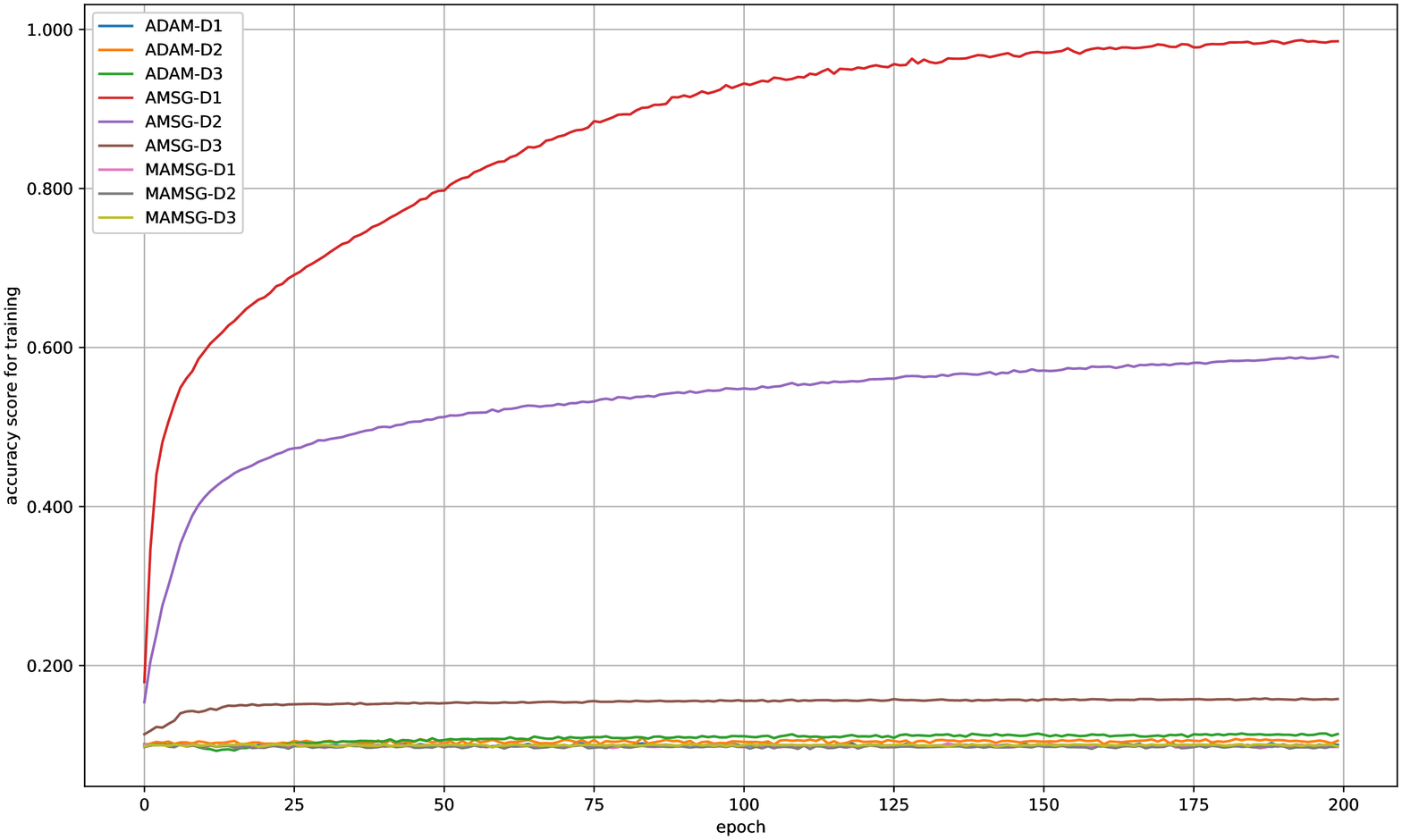}
\label{fig:2_d_a}
}
\hfil
\subfloat[]{
\includegraphics[width=9truecm]{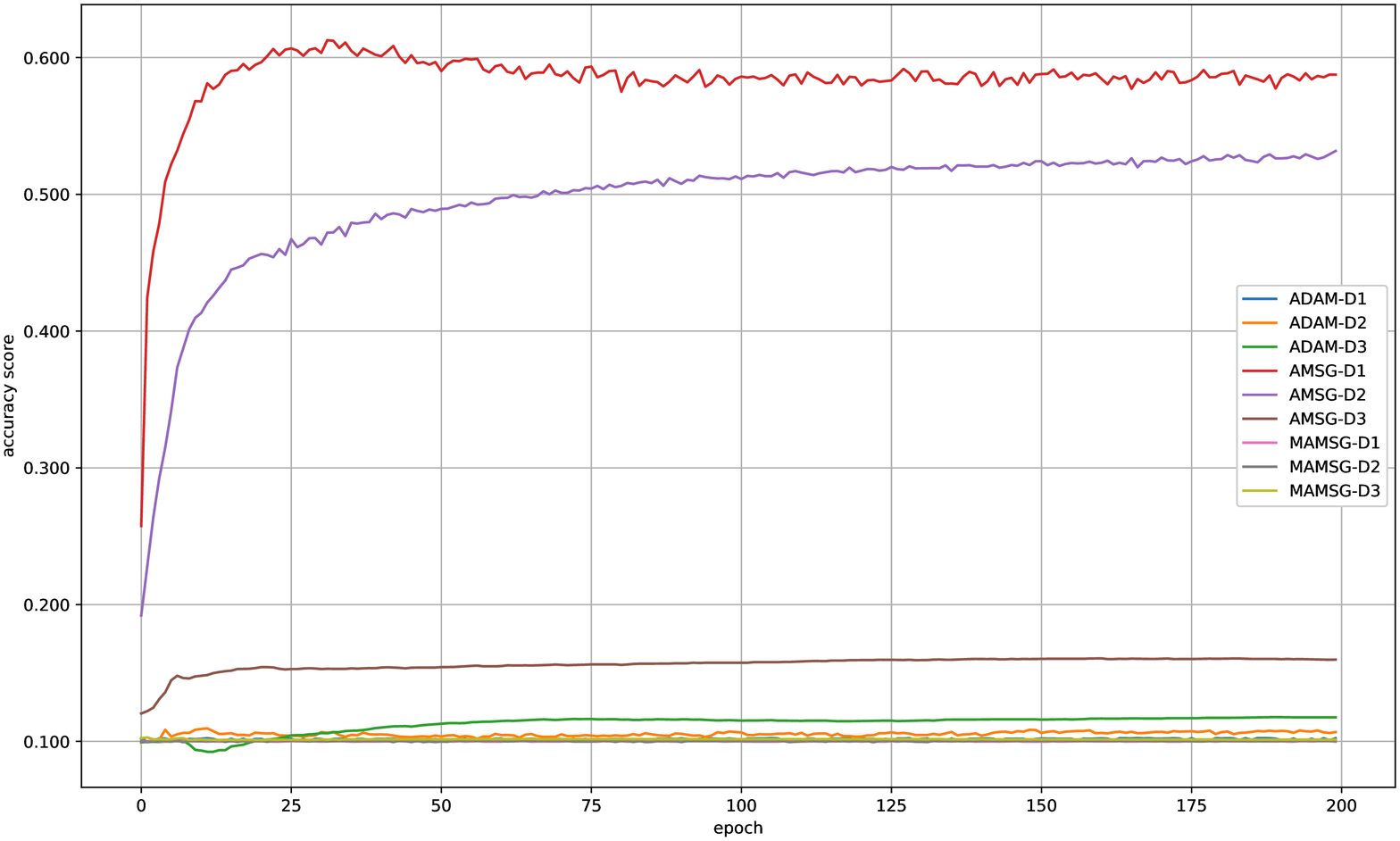}
\label{fig:2_d_a_t}
}
\caption{(a) Training loss function value, (b) training classification accuracy score, and (c) test classification accuracy score for Algorithm \ref{algo:1} with diminishing sub-learning rates versus number of epochs on the CIFAR-10 dataset}
\label{fig:2_d}
\end{figure*}

\section{Conclusion}
\label{sec:6}
This paper proposed an adaptive learning rate optimization algorithm for solving {stationary point problems associated with nonconvex} stochastic optimization problems in deep learning and performed convergence and convergence rate analyses for constant sub-learning rates and diminishing sub-learning rates. {In particular,} the analyses show that the proposed algorithm with constant sub-learning rates can solve the problem.
Numerical results showed that the proposed algorithm can {be applied to stochastic optimization in text and image classification tasks}, while Adam and AMSGrad with diminishing sub-learning rates cannot do so.  It also had higher classification accuracy compared with the existing algorithms. In particular, the results showed that the proposed algorithm using constant sub-learning rates is well suited to training neural networks.

\section{Proofs of Theorems \ref{thm:1} and \ref{thm:2} {and Propositions \ref{prop:0} and \ref{prop:1}}}
\label{appen:1}
The history of the process $\bm{\xi}_0,\bm{\xi}_1,\ldots$ up to time $n$ is denoted by $\bm{\xi}_{[n]} = (\bm{\xi}_0,\bm{\xi}_1,\ldots,\bm{\xi}_n)$. For a random process $\bm{\xi}_0, \bm{\xi}_1, \ldots$, let $\mathbb{E}[Y|\bm{\xi}_{[n]}]$ denote the conditional expectation of $Y$ given $\bm{\xi}_{[n]} = (\bm{\xi}_0,\bm{\xi}_1,\ldots,\bm{\xi}_n)$. Let $(\bm{x}_n)_{n\in\mathbb{N}}$, $(\bm{m}_n)_{n\in\mathbb{N}}$, $(\bm{\mathsf{d}}_n)_{n\in\mathbb{N}}$, and $(\mathsf{H}_n)_{n\in\mathbb{N}}$ be the sequences generated by Algorithm \ref{algo:1}. First, we prove a lemma.

\begin{lem}\label{lem:1}
Suppose that (A1)--(A2) and (C1)--(C2) hold. Then, for all $\bm{x} \in X$ and all $n\in\mathbb{N}$,
\begin{align*}
\mathbb{E}\left[\left\| \bm{x}_{n+1} - \bm{x} \right\|_{\mathsf{H}_n}^2\right]
&\leq
\mathbb{E}\left[\left\| \bm{x}_{n} - \bm{x} \right\|_{\mathsf{H}_n}^2\right]
+ 2 \alpha_n 
 \bigg\{
 \frac{1 - \beta_n}{1 - {\gamma}^{n+1}} 
 {\mathbb{E} \left[\left\langle \bm{x} - \bm{x}_n, \nabla f (\bm{x}_n) \right\rangle \right]}\\
&\quad + \frac{\beta_n}{1 - {\gamma}^{n+1}} \mathbb{E} \left[ \left\langle \bm{x} - \bm{x}_n,\bm{m}_{n-1} \right\rangle \right]
 \bigg\}
 + \alpha_n^2 \mathbb{E} \left[ \left\|\bm{\mathsf{d}}_n \right\|_{\mathsf{H}_n}^2 \right]. 
\end{align*}
\end{lem}

{\em Proof:} Choose $\bm{x}\in X$ and $n\in \mathbb{N}$. From the definition of $\bm{x}_{n+1}$ and the nonexpansivity of $P_{X,\mathsf{H}_n}$ {(i.e., $\|P_{X,\mathsf{H}_n}(\bm{x}) - P_{X,\mathsf{H}_n}(\bm{y}) \|_{\mathsf{H}_n} \leq \|\bm{x} - \bm{y} \|_{\mathsf{H}_n}$ ($\bm{x},\bm{y} \in \mathbb{R}^d$))}, we have, almost surely,
\begin{align*}
\left\|\bm{x}_{n+1} - \bm{x} \right\|_{\mathsf{H}_n}^2
&\leq
\left\| (\bm{x}_n - \bm{x}) + \alpha_n \bm{\mathsf{d}}_n \right\|_{\mathsf{H}_n}^2\\
&=
\left\| \bm{x}_n - \bm{x} \right\|_{\mathsf{H}_n}^2 
+ 2 \alpha_n \left\langle \bm{x}_n - \bm{x}, \bm{\mathsf{d}}_n \right\rangle_{\mathsf{H}_n}
+ \alpha_n^2 \left\|\bm{\mathsf{d}}_n \right\|_{\mathsf{H}_n}^2. 
\end{align*}
Moreover, the definitions of $\bm{\mathsf{d}}_n$, $\bm{m}_n$, and $\hat{\bm{m}}_n$ ensure that
\begin{align*}
\left\langle \bm{x}_n - \bm{x}, \bm{\mathsf{d}}_n \right\rangle_{\mathsf{H}_n}=
\frac{1}{{\tilde{\gamma}}_n}\left\langle \bm{x} - \bm{x}_n, \bm{m}_n \right\rangle
=
\frac{\beta_n}{{\tilde{\gamma}}_n} \left\langle \bm{x} - \bm{x}_n, \bm{m}_{n-1} \right\rangle 
+
\frac{1-\beta_n}{{\tilde{\gamma}}_n} \left\langle \bm{x} - \bm{x}_n, \mathsf{G}(\bm{x}_n,\bm{\xi}_n) \right\rangle,
\end{align*}
where ${\tilde{\gamma}}_n := 1 - {\gamma}^{n+1}$. Hence, almost surely,
\begin{align}
\left\|\bm{x}_{n+1} - \bm{x} \right\|_{\mathsf{H}_n}^2
&\leq
\left\| \bm{x}_n -\bm{x} \right\|_{\mathsf{H}_n}^2
+ 
2 \alpha_n \bigg\{
\frac{\beta_n}{{\tilde{\gamma}}_n} \left\langle \bm{x} - \bm{x}_n, \bm{m}_{n-1} \right\rangle 
+
\frac{1-\beta_n}{{\tilde{\gamma}}_n} \left\langle \bm{x} - \bm{x}_n, \mathsf{G}(\bm{x}_n,\bm{\xi}_n) \right\rangle
\bigg\}\nonumber \\
&\quad + \alpha_n^2 \left\| \bm{\mathsf{d}}_n \right\|_{\mathsf{H}_n}^2.
\label{ineq:004}
\end{align}
The condition $\bm{x}_n = \bm{x}_n(\bm{\xi}_{[n-1]})$ $(n\in \mathbb{N})$, {(C1), and (C2)} guarantee that
\begin{align*}
\mathbb{E} \left[\left\langle \bm{x} - \bm{x}_n, \mathsf{G}(\bm{x}_n,\bm{\xi}_n) \right\rangle \right]
&=
\mathbb{E} 
 \left[ 
 \mathbb{E} \left[\left\langle \bm{x} - \bm{x}_n, \mathsf{G}(\bm{x}_n,\bm{\xi}_n) \right\rangle | \bm{\xi}_{[n-1]} \right]
 \right]\\
&= 
\mathbb{E} 
 \left[ 
 \left\langle \bm{x} - \bm{x}_n, \mathbb{E} \left[ \mathsf{G}(\bm{x}_n,\bm{\xi}_n) | \bm{\xi}_{[n-1]} \right] \right\rangle
 \right]\\
&=
\mathbb{E} \left[\left\langle \bm{x} - \bm{x}_n, {\nabla f} (\bm{x}_n) \right\rangle \right].
\end{align*}
Therefore, the lemma follows by taking the expectation of \eqref{ineq:004}.
\qed

\begin{lem}\label{lem:2}
If (C3) holds, then, for all $n\in\mathbb{N}$, $\mathbb{E}[\|\bm{m}_n\|^2] \leq \tilde{M}^2 := \max \{ \|\bm{m}_{-1}\|^2, M^2 \}$. Moreover, if {(A3)} holds, then, for all $n\in\mathbb{N}$, $\mathbb{E}[\|\bm{\mathsf{d}}_n\|_{\mathsf{H}_n}^2] \leq \tilde{B}^2 \tilde{M}^2/(1-{\gamma})^2$, where $\tilde{B} := \sup\{ {\max_{i=1,2,\ldots,d} h_{n,i}^{-1/2}} \colon n\in\mathbb{N}\} < + \infty$. 
\end{lem}

{\em Proof:} The convexity of $\|\cdot\|^2$, together with the definition of $\bm{m}_n$ and (C3), guarantees that, for all $n\in\mathbb{N}$,
\begin{align*}
\mathbb{E}\left[ \left\|\bm{m}_n\right\|^2 \right]
&\leq 
\beta_n \mathbb{E}\left[ \left\|\bm{m}_{n-1} \right\|^2 \right] + (1-\beta_n) M^2.
\end{align*}
Induction thus ensures that, for all $n\in\mathbb{N}$,
\begin{align}\label{induction}
\mathbb{E} \left[ \left\|\bm{m}_n \right\|^2 \right] \leq \tilde{M}^2 := \max \left\{ \left\|\bm{m}_{-1}\right\|^2,M^2 \right\} < + \infty.
\end{align}
Given $n\in\mathbb{N}$, $\mathsf{H}_n \succ O$ ensures that there exists a unique matrix $\overline{\mathsf{H}}_n \succ O$ such that $\mathsf{H}_n = \overline{\mathsf{H}}_n^2$ \cite[Theorem 7.2.6]{horn}. 
From $\|\bm{x}\|_{\mathsf{H}_n}^2 = \| \overline{\mathsf{H}}_n \bm{x} \|^2$ for all $\bm{x}\in\mathbb{R}^d$ and the definitions of $\bm{\mathsf{d}}_n$ and $\hat{\bm{m}}_n$, we have that, for all $n\in\mathbb{N}$, 
\begin{align*}
\mathbb{E} \left[ \left\| \bm{\mathsf{d}}_n \right\|_{\mathsf{H}_n}^2 \right]
= 
\mathbb{E} \left[ \left\| \overline{\mathsf{H}}_n^{-1} \mathsf{H}_n\bm{\mathsf{d}}_n \right\|^2 \right]
\leq 
\frac{1}{{\tilde{\gamma}}_n^2}\mathbb{E} \left[ \left\| \overline{\mathsf{H}}_n^{-1} \right\|^2 \left\|\bm{m}_n \right\|^2 \right],
\end{align*}
where ${\tilde{\gamma}}_n := 1 - {\gamma}^{n+1} \geq 1 - {\gamma}$ and $\| \overline{\mathsf{H}}_n^{-1} \| = \| \mathsf{diag}(h_{n,i}^{-1/2}) \| = {\max_{i=1,2,\ldots,d} h_{n,i}^{-1/2}}$ ($n\in \mathbb{N}$). From \eqref{induction} and $\tilde{B} := \sup \{ {\max_{i=1,2,\ldots,d} h_{n,i}^{-1/2}} \colon n\in\mathbb{N} \} \leq \max_{i=1,2,\ldots,d} {h_{0,i}^{-1/2}} < + \infty$ (by {(A3)}), we have that, for all $n\in \mathbb{N}$, $\mathbb{E} [\| \bm{\mathsf{d}}_n \|_{\mathsf{H}_n}^2 ] \leq \tilde{B}^2 \tilde{M}^2/(1-{\gamma})^2$, which completes the proof.
\qed

The convergence rate analysis of Algorithm \ref{algo:1} is as follows.

\begin{thm}\label{thm:3}
Suppose that (A1)--(A5) and (C1)--(C3) hold and $(\gamma_n)_{n\in \mathbb{N}}$ defined by $\gamma_n := \alpha_n (1-\beta_n)/(1 - {\gamma}^{n+1})$ and $(\beta_n)_{n\in \mathbb{N}}$ satisfy $\gamma_{n+1} \leq \gamma_n$ ($n\in\mathbb{N}$) and $\limsup_{n\to + \infty} \beta_n < 1$. {Let $V_n (\bm{x}) = V_n := \mathbb{E}\left[ \langle \bm{x}_n - \bm{x}, \nabla f (\bm{x}_n) \rangle \right]$ for all $\bm{x} \in X$ and all $n\in\mathbb{N}$.} Then, {for all $\bm{x} \in X$ and all $n\geq 1$,}
\begin{align*}
{\frac{1}{n} \sum_{k=1}^n V_k}
\leq
\frac{D \sum_{i=1}^d B_i}{2 \tilde{b} n \alpha_n}
+
\frac{\tilde{B}^2 \tilde{M}^2}{2 \tilde{b} {\tilde{\gamma}}^2 n}\sum_{k=1}^n \alpha_k
+
\frac{\tilde{M}\sqrt{Dd}}{\tilde{b}n} \sum_{k=1}^n \beta_k,
\end{align*}
where $(\beta_n)_{n\in\mathbb{N}} \subset (0,b] \subset (0,1)$, $\tilde{b} := 1 - b$, {$\tilde{\gamma} := 1-\gamma$}, $\tilde{M}$ and $\tilde{B}$ are defined as in Lemma \ref{lem:2}, and $D$ and $B_i$ are defined as in the Assumption.
\end{thm}

{\em Proof:} Let $\bm{x}\in X$ be fixed arbitrarily. Lemma \ref{lem:1} guarantees that, for all $k\in\mathbb{N}$,
\begin{align*}
{V_k}
&\leq
\frac{1}{2\gamma_k}
\left\{ 
\mathbb{E}\left[\left\| \bm{x}_{k} - \bm{x} \right\|_{\mathsf{H}_k}^2\right]
-
\mathbb{E}\left[\left\| \bm{x}_{k+1} - \bm{x} \right\|_{\mathsf{H}_k}^2\right]
\right\}\\
&\quad + \frac{\beta_k}{1 - \beta_k} \mathbb{E} \left[ \left\langle \bm{x} - \bm{x}_k, \bm{m}_{k-1} \right\rangle \right]
+ \frac{\alpha_k {\tilde{\gamma}}_k}{2 (1-\beta_k)} \mathbb{E} \left[ \left\|\bm{\mathsf{d}}_k \right\|_{\mathsf{H}_k}^2 \right],
\end{align*}
where ${\tilde{\gamma}}_n := 1 - {\gamma}^{n+1} \leq 1$ ($n\in\mathbb{N}$). The condition $\limsup_{n\to + \infty} \beta_n < 1$ ensures the existence of $b > 0$ such that, for all $n\in \mathbb{N}$, $\beta_n \leq b < 1$. Let $\tilde{b} := 1 -b$. Then, for all $n \geq 1$,
\begin{align}\label{key}
{\sum_{k=1}^n V_k}\nonumber
&\leq
\frac{1}{2} \underbrace{\sum_{k=1}^n \frac{1}{\gamma_k}
\left\{ 
\mathbb{E}\left[\left\| \bm{x}_{k} - \bm{x} \right\|_{\mathsf{H}_k}^2\right]
-
\mathbb{E}\left[\left\| \bm{x}_{k+1} - \bm{x} \right\|_{\mathsf{H}_k}^2\right]
\right\}}_{\Gamma_n}\nonumber\\
&\quad + \underbrace{
\sum_{k=1}^n \frac{\beta_k}{1-\beta_k} \mathbb{E} \left[ \left\langle \bm{x} - \bm{x}_k,\bm{m}_{k-1} \right\rangle \right]}_{B_n}\nonumber\\
&\quad + \frac{1}{2 \tilde{b}} \underbrace{\sum_{k=1}^n \alpha_k \mathbb{E} \left[ \left\|\bm{\mathsf{d}}_k \right\|_{\mathsf{H}_k}^2 \right]}_{A_n}.
\end{align} 
From the definition of $\Gamma_n$ and $\mathbb{E} [ \| \bm{x}_{n+1} - \bm{x} \|_{\mathsf{H}_{n}}^2]/\gamma_n \geq 0$, 
\begin{align}\label{LAM}
\Gamma_n
&\leq
\frac{\mathbb{E}\left[\left\| \bm{x}_{1} - \bm{x} \right\|_{\mathsf{H}_{1}}^2\right]}{\gamma_1}
+
\underbrace{
\sum_{k=2}^n \left\{
\frac{\mathbb{E}\left[\left\| \bm{x}_{k} - \bm{x} \right\|_{\mathsf{H}_{k}}^2\right]}{\gamma_k}
-
\frac{\mathbb{E}\left[\left\| \bm{x}_{k} - \bm{x} \right\|_{\mathsf{H}_{k-1}}^2\right]}{\gamma_{k-1}} 
\right\}
}_{\tilde{\Gamma}_n}.
\end{align}
Since $\overline{\mathsf{H}}_k \succ O$ exists such that $\mathsf{H}_k = \overline{\mathsf{H}}_k^2$, we have $\|\bm{x}\|_{\mathsf{H}_k}^2 = \| \overline{\mathsf{H}}_k \bm{x} \|^2$ for all $\bm{x}\in\mathbb{R}^d$. Accordingly, we have 
\begin{align*}
\tilde{\Gamma}_n 
=
\mathbb{E} \left[ 
\sum_{k=2}^n 
\left\{
\frac{\left\| \overline{\mathsf{H}}_{k} (\bm{x}_{k} - \bm{x}) \right\|^2}{\gamma_k}
-
\frac{\left\| \overline{\mathsf{H}}_{k-1} (\bm{x}_{k} - \bm{x}) \right\|^2}{\gamma_{k-1}}
\right\}
\right].
\end{align*}
The Assumption ensures we can express $\mathsf{H}_k$ as $\mathsf{H}_k = \mathsf{diag}(h_{k,i})$, where $h_{k,i} > 0$ $(k\in\mathbb{N}, i=1,2,\ldots,d)$. Hence, for all $k\in\mathbb{N}$ and all $\bm{x} := (x_i) \in \mathbb{R}^d$,
\begin{align}\label{HK}
\overline{\mathsf{H}}_{k} = \mathsf{diag}\left(\sqrt{h_{k,i}}\right) 
\text{ and }
\left\| \overline{\mathsf{H}}_{k} \bm{x} \right\|^2
=
\sum_{i=1}^d h_{k,i} x_i^2.
\end{align}
Hence, for all $n\geq 2$,
\begin{align*}
\tilde{\Gamma}_n 
= 
\mathbb{E} \left[ 
\sum_{k=2}^n
\sum_{i=1}^d 
\left(
\frac{h_{k,i}}{\gamma_k}
-
\frac{h_{k-1,i}}{\gamma_{k-1}}
\right)
(x_{k,i} - x_i)^2
\right].
\end{align*}
From $\gamma_k \leq \gamma_{k-1}$ $(k\geq 1)$ and (A3), we have $h_{k,i}/\gamma_k - h_{k-1,i}/\gamma_{k-1} \geq 0$ $(k \geq 1, i=1,2,\ldots,d)$. Moreover, from (A5), $\max_{i=1,2,\ldots,d} \sup \{ (x_{n,i} - x_i)^2 \colon n\in\mathbb{N} \} \leq D < + \infty$. Accordingly, for all $n\geq 2$,
\begin{align*}
\tilde{\Gamma}_n
\leq
D
\mathbb{E} \left[ 
\sum_{k=2}^n
\sum_{i=1}^d 
\left(
\frac{h_{k,i}}{\gamma_k}
-
\frac{h_{k-1,i}}{\gamma_{k-1}}
\right)
\right]
= 
D
\mathbb{E} \left[ 
\sum_{i=1}^d
\left(
\frac{h_{n,i}}{\gamma_n}
-
\frac{h_{1,i}}{\gamma_{1}}
\right)
\right].
\end{align*}
Therefore, \eqref{LAM}, $\mathbb{E} [\| \bm{x}_{1} - \bm{x}\|_{\mathsf{H}_{1}}^2]/\gamma_1 \leq D \mathbb{E} [ \sum_{i=1}^d h_{1,i}/\gamma_1]$, and (A4) imply, for all $n\in\mathbb{N}$,
\begin{align*}
\Gamma_n 
&\leq
D \mathbb{E} \left[ \sum_{i=1}^d \frac{h_{1,i}}{\gamma_1} \right]
+
D
\mathbb{E} \left[
\sum_{i=1}^d 
\left(
\frac{h_{n,i}}{\gamma_n}
-
\frac{h_{1,i}}{\gamma_{1}}
\right)
\right]\\
&=
\frac{D}{\gamma_n}
\mathbb{E} \left[
\sum_{i=1}^d 
h_{n,i}
\right]
\leq 
\frac{D}{\gamma_n}
\sum_{i=1}^d 
B_{i},
\end{align*}
which, together with $\gamma_n := \alpha_n(1-\beta_n)/(1 - {\gamma}^{n+1})$ and $\tilde{b} := 1 -b$, implies 
\begin{align}\label{L} 
\Gamma_n 
\leq 
\frac{D \sum_{i=1}^d 
B_{i}}{\tilde{b} \alpha_n}.
\end{align}
The Cauchy-Schwarz inequality, together with $\max_{i=1,2,\ldots,d} \sup \{ (x_{n,i} - x_i)^2 \colon n\in\mathbb{N} \} \leq D < + \infty$ (by (A5)) and $\mathbb{E}[\|\bm{m}_n\|] \leq \tilde{M}$ $(n\in\mathbb{N})$ (by Lemma \ref{lem:2}), guarantees that, for all $n\in\mathbb{N}$,
\begin{align}\label{B}
\begin{split}
B_n 
\leq 
\frac{\sqrt{Dd}}{\tilde{b}} \sum_{k=1}^n \beta_k \mathbb{E} \left[
\left\|\bm{m}_{k-1} \right\|
\right]
\leq
\frac{\tilde{M}\sqrt{Dd}}{\tilde{b}} \sum_{k=1}^n \beta_k.
\end{split}
\end{align}
Since $\mathbb{E}[ \|\bm{\mathsf{d}}_n \|_{\mathsf{H}_n}^2 ] \leq \tilde{B}^2 \tilde{M}^2/(1 - {\gamma})^2$ $(n\in\mathbb{N})$ holds (by Lemma \ref{lem:2}), we have, for all $n\in\mathbb{N}$,
\begin{align}\label{D}
A_n := \sum_{k=1}^n \alpha_k \mathbb{E} \left[ \left\|\bm{\mathsf{d}}_k \right\|_{\mathsf{H}_k}^2 \right] 
\leq \frac{\tilde{B}^2 \tilde{M}^2}{(1-{\gamma})^2} \sum_{k=1}^n \alpha_k.
\end{align}
Therefore, \eqref{key}, \eqref{L}, \eqref{B}, and \eqref{D} leads to {the assertion in Theorem \ref{thm:3}}. This completes the proof.
\qed

{Lemmas \ref{lem:1} and \ref{lem:2} lead to Theorem \ref{thm:1}.}

{\em {Proof of Theorem \ref{thm:1}:}} { Let $\bm{x} \in X$, $\alpha_n := \alpha \in (0,1)$, and $\beta_n := \beta = b \in (0,1)$. We show that, for all $\epsilon > 0$,
\begin{align}\label{Ineq:1}
\liminf_{n\to +\infty} 
V_n
\leq 
\frac{\tilde{B}^2 \tilde{M}^2}{2\tilde{b} \tilde{\gamma}^2} \alpha
+
\frac{\tilde{M}\sqrt{Dd}}{\tilde{b}\tilde{\gamma}}\beta
+
\frac{Dd\epsilon}{2 \tilde{b}}
+ 
\epsilon.
\end{align}
If \eqref{Ineq:1} does not hold for all $\epsilon > 0$, then there exists $\epsilon_0 > 0$ such that 
\begin{align}\label{INF}
\liminf_{n\to +\infty} 
V_n
> 
\frac{\tilde{B}^2 \tilde{M}^2}{2\tilde{b}\tilde{\gamma}^2} \alpha
+
\frac{\tilde{M}\sqrt{Dd}}{\tilde{b}\tilde{\gamma}}\beta
+
\frac{Dd\epsilon_0}{2 \tilde{b}}
+ 
\epsilon_0.
\end{align}
Assumptions (A3) and (A4) ensure that there exists $n_0\in\mathbb{N}$ such that, for all $n \in \mathbb{N}$, $n\geq n_0$ implies that 
\begin{align}\label{I}
\mathbb{E}\left[\sum_{i=1}^d (h_{n+1,i} - h_{n,i}) \right] 
\leq \frac{d \alpha \epsilon_0}{2}.
\end{align}
From \eqref{HK}, (A3), (A5), and \eqref{I}, for all $n \geq n_0$,
\begin{align}\label{IV}
X_{n+1}
-
\mathbb{E}\left[ \left\| \bm{x}_{n+1} - \bm{x} \right\|_{\mathsf{H}_{n}}^2 \right]
= 
\mathbb{E} \left[
\sum_{i=1}^d (h_{n+1,i} - h_{n,i})(x_{n+1,i} - x_i)^2
\right]
\leq
\frac{D d \alpha \epsilon_0}{2},
\end{align}
where, for all $n\in\mathbb{N}$, $X_n := \mathbb{E} [ \| \bm{x}_{n} - \bm{x} \|_{\mathsf{H}_{n}}^2 ] \leq D \sum_{i=1}^d B_i < + \infty$ from (A4) and (A5). Moreover, from $\gamma \in [0,1)$, there exists $n_1 \in \mathbb{N}$ such that, for all $n \in \mathbb{N}$, $n \geq n_1$ implies that 
\begin{align}\label{I_I}
X_{n+1} \gamma^{n+1} 
\leq 
\frac{D d \alpha \epsilon_0}{2}.
\end{align}
The definition of the limit inferior of $(V_n)_{n\in\mathbb{N}}$ guarantees that there exists $n_2 \in \mathbb{N}$ such that, for all $n \geq n_2$,
\begin{align*}
\liminf_{n\to +\infty} 
V_n
-\frac{1}{2} \epsilon_0 
\leq 
V_n,
\end{align*}
which, together with \eqref{INF}, implies that, for all $n \geq n_1$,
\begin{align}\label{II}
V_n
> 
\frac{\tilde{B}^2 \tilde{M}^2}{2\tilde{b}\tilde{\gamma}^2} \alpha
+
\frac{\tilde{M}\sqrt{Dd}}{\tilde{b}\tilde{\gamma}}\beta
+
\frac{Dd\epsilon_0}{2 \tilde{b}}
+ 
\frac{1}{2} \epsilon_0.
\end{align}
Lemmas \ref{lem:1} and \ref{lem:2} and \eqref{IV} thus lead to the finding that, for all $n \geq n_3 := \max \{n_0,n_1,n_2\}$,
\begin{align*}
X_{n+1} 
\leq 
X_n + \frac{D d \alpha \epsilon_0}{2}
- \frac{2 \alpha \tilde{b}}{1-\gamma^{n+1}} V_n
+ \frac{2\tilde{M}\sqrt{Dd}}{\tilde{\gamma}}\alpha \beta
 + \frac{\tilde{B}^2 \tilde{M}^2}{\tilde{\gamma}^2}\alpha^2,
\end{align*}
where $\tilde{b} := 1 -b$ and $\tilde{\gamma} := 1 - \gamma$. Hence, from \eqref{I_I}, $1 - \gamma^{n+1} \leq 1$, and $(X_{n+1} - X_n) \gamma^{n+1} \leq X_{n+1} \gamma^{n+1}$ ($n\in\mathbb{N}$), we have that, for all $n \geq n_3$,
\begin{align}\label{In:1}
\begin{split}
X_{n+1} 
&\leq 
X_n + \frac{D d \alpha \epsilon_0}{2}
- 2 \alpha \tilde{b} V_n
+ \frac{2\tilde{M}\sqrt{Dd}}{\tilde{\gamma}}\alpha \beta
 + \frac{\tilde{B}^2 \tilde{M}^2}{\tilde{\gamma}^2}\alpha^2
+ X_{n+1} \gamma^{n+1}\\
&\leq 
X_n + D d \alpha \epsilon_0
- 2 \alpha \tilde{b} V_n
+ \frac{2\tilde{M}\sqrt{Dd}}{\tilde{\gamma}}\alpha \beta
 + \frac{\tilde{B}^2 \tilde{M}^2}{\tilde{\gamma}^2}\alpha^2.
\end{split}
\end{align}
Therefore, \eqref{II} ensures that, for all $n \geq n_3$,
\begin{align*}
X_{n+1} 
&<
X_n + D d \alpha \epsilon_0
-2 \alpha \tilde{b} 
\bigg\{
\frac{\tilde{B}^2 \tilde{M}^2}{2\tilde{b} \tilde{\gamma}^2} \alpha
+
\frac{\tilde{M}\sqrt{Dd}}{\tilde{b}\tilde{\gamma}}\beta
+
\frac{Dd\epsilon_0}{2 \tilde{b}}
+ 
\frac{1}{2} \epsilon_0
\bigg\}
+ \frac{2\tilde{M}\sqrt{Dd}}{\tilde{\gamma}}\alpha \beta
+ \frac{\tilde{B}^2 \tilde{M}^2}{\tilde{\gamma}^2}\alpha^2\\
&=
X_n - \alpha \tilde{b} \epsilon_0\\
&< 
X_{n_3} - \alpha \tilde{b} \epsilon_0 (n +1 - n_3).
\end{align*}
Since the right-hand side of the above inequality approaches minus infinity when $n$ diverges, we have a contradiction. Hence, \eqref{Ineq:1} holds for all $\epsilon > 0$. From the arbitrary condition of $\epsilon$, we have that
\begin{align*}
\liminf_{n\to +\infty} 
V_n 
\leq 
\frac{\tilde{B}^2 \tilde{M}^2}{2\tilde{b}\tilde{\gamma}^2} \alpha
+
\frac{\tilde{M}\sqrt{Dd}}{\tilde{b}\tilde{\gamma}}\beta,
\end{align*}
}
which completes the proof.
\qed

{Lemmas \ref{lem:1} and \ref{lem:2} and Theorem \ref{thm:3} lead to Theorem \ref{thm:2}.}

{\em {Proof of Theorem \ref{thm:2}:}} { Lemmas \ref{lem:1} and \ref{lem:2}, together with a discussion similar to the one for obtaining \eqref{In:1}, ensure that, for all $k\in\mathbb{N}$,
\begin{align*}
X_{k+1} 
&\leq 
X_k
+ D \mathbb{E}\left[ \sum_{i=1}^d (h_{k+1,i} - h_{k,i}) \right]
- 2 \alpha_k (1-\beta_k) V_k\\
&\quad 
+ \frac{2 \tilde{M} \sqrt{Dd}}{\tilde{\gamma}}\alpha_k \beta_k 
+ \frac{\tilde{B}^2 \tilde{M}^2}{\tilde{\gamma}^2}\alpha_k^2
+ D \sum_{i=1}^d B_i \gamma^{k+1},
\end{align*} 
which implies that
\begin{align*}
2 \alpha_k V_k
&\leq 
X_k - X_{k+1}
+ D \mathbb{E}\left[ \sum_{i=1}^d (h_{k+1,i} - h_{k,i}) \right]
+ \frac{\tilde{B}^2 \tilde{M}^2}{\tilde{\gamma}^2}\alpha_k^2\\
&\quad + 2 
\left(
\frac{\tilde{M} \sqrt{Dd}}{\tilde{\gamma}}
+ F
\right)
\alpha_k \beta_k
+ D \sum_{i=1}^d B_i \gamma^{k+1}, 
\end{align*}
where $F := \sup \{ |V_n| \colon n\in\mathbb{N} \} < + \infty$ holds from (A2) and (A5). Summing up the above inequality from $k=0$ to $k=n$ ensures that
\begin{align*}
2 \sum_{k=0}^n \alpha_k V_k
&\leq 
X_0
+ D \mathbb{E}\left[ \sum_{i=1}^d (h_{n+1,i} - h_{0,i}) \right]
+ \frac{\tilde{B}^2 \tilde{M}^2}{\tilde{\gamma}^2} \sum_{k=0}^n \alpha_k^2\\
&\quad + 2 
\left(
\frac{\tilde{M} \sqrt{Dd}}{\tilde{\gamma}}
+ F
\right)
\sum_{k=0}^n \alpha_k \beta_k
+ D \hat{B} \sum_{k=0}^n \gamma^{k+1},
\end{align*}
where $\hat{B} := \sum_{i=1}^d B_i$. Let $(\alpha_n)_{n\in\mathbb{N}}$ and $(\beta_n)_{n\in\mathbb{N}}$ satisfy $\sum_{n=0}^{+\infty} \alpha_n = + \infty$, $\sum_{n=0}^{+\infty} \alpha_n^2 < + \infty$, and $\sum_{n=0}^{+\infty} \alpha_n \beta_n < + \infty$. From (A4) and $\gamma \in [0,1)$, we have that 
\begin{align}\label{sum}
\sum_{k=0}^{+\infty} \alpha_k V_k
< + \infty.
\end{align}
We prove that $\liminf_{n \to +\infty} V_n \leq 0$. If $\liminf_{n \to +\infty} V_n \leq 0$ does not hold, then there exist $\zeta > 0$ and $m_0\in\mathbb{N}$ such that, for all $n \geq m_0$, $V_n \geq \zeta$. Accordingly, \eqref{sum} and $\sum_{n=0}^{+\infty} \alpha_n = + \infty$ guarantee that
\begin{align*}
+ \infty = 
\zeta \sum_{k=m_0}^{+\infty} \alpha_k \leq \sum_{k=m_0}^{+\infty} \alpha_k V_k
< + \infty,
\end{align*}
which is a contradiction. Hence, $\liminf_{n \to +\infty} V_n \leq 0$ holds. }

Let $\alpha_n := 1/n^\eta$ ($\eta \in [1/2,1)$) and let $\beta_n := \lambda^n$ ($\lambda \in (0,1)$). Then, $\gamma_{n+1} \leq \gamma_n$ ($n\in\mathbb{N}$) and $\limsup_{n\to + \infty} \beta_n < 1$. 
We have that 
\begin{align*}
\lim_{n \to + \infty} \frac{1}{n \alpha_n} = \lim_{n \to + \infty} \frac{1}{n^{1-\eta}} = 0
\end{align*}
and 
\begin{align*}
\frac{1}{n} \sum_{k=1}^n \alpha_k
\leq 
\frac{1}{n} \left\{ 1 + \int_1^n \frac{\mathrm{d}t}{t^\eta}  \right\}
= \frac{1}{n} \left\{ \frac{n^{1-\eta}}{1 - \eta} - \frac{\eta}{1-\eta}  \right\}
\leq \frac{1}{(1-\eta) n^\eta}
\leq \frac{1}{(1-\eta) n^{1-\eta}}.
\end{align*}
Therefore, Theorem \ref{thm:3} ensures that $(1/n) \sum_{k=1}^n V_k \leq \mathcal{O} ( 1/n^{1 -\eta} )$, which completes the proof.
\qed

{\em {Proof of Proposition \ref{prop:0}:}} { Since $F(\cdot,\bm{\xi})$ is convex for almost every $\bm{\xi} \in \Xi$, we have that, for all $n\in\mathbb{N}$, $\mathbb{E} [f(\tilde{\bm{x}}_n) - f^\star ] \leq (1/n) \sum_{k=1}^n \mathbb{E} [f(\bm{x}_k) - f^\star] \leq (1/n) \sum_{k=1}^n V_k$, which, together with Theorem \ref{thm:1}, leads to Proposition \ref{prop:0}.
\qed
}

{\em {Proof of Proposition \ref{prop:1}:}} { Theorem \ref{thm:2} and the proof of Proposition \ref{prop:0} lead to the finding that $\lim_{n \to + \infty} \mathbb{E} [ f ( \tilde{\bm{x}}_n ) - f^\star] = 0$. Let $\hat{\bm{x}} \in X$ be an arbitrary accumulation point of $(\tilde{\bm{x}}_n)_{n\in\mathbb{N}} \subset X$. Since there exists $(\tilde{\bm{x}}_{n_i})_{i\in\mathbb{N}} \subset (\tilde{\bm{x}}_n)_{n\in\mathbb{N}}$ such that $(\tilde{\bm{x}}_{n_i})_{i\in\mathbb{N}}$ converges almost surely to $\hat{\bm{x}}$, the continuity of $f$ and $\lim_{n \to + \infty} \mathbb{E} [ f ( \tilde{\bm{x}}_n ) - f^\star] = 0$ imply that $\mathbb{E} \left[ f ( \hat{\bm{x}}) - f^\star \right]= 0$, and hence, $\hat{\bm{x}} \in X^\star$.
\qed
}

\section*{Acknowledgment} 
I thank Hiroyuki Sakai for his input on the numerical examples.

\end{document}